\documentclass[12pt,a4paper]{article}

\usepackage[utf8]{inputenc}  
\usepackage[T1]{fontenc}
\usepackage{lmodern}

\usepackage{amsmath,amssymb,amsfonts,amsthm}
\usepackage{graphicx}
\usepackage{bm}
\usepackage{xcolor}
\usepackage{hyperref}
\usepackage[numbers]{natbib}
\usepackage{mathabx}

\usepackage[ruled,vlined]{algorithm2e}

\newtheorem{remark}{Remark}
\newtheorem{theorem}{Theorem}
\newtheorem{defi}{Definition}
\newtheorem{proposition}{Proposition}
\newtheorem{property}{Property}
\newtheorem{lemma}{Lemma}
\newtheorem{corollary}{Corollary}

\DeclareMathOperator{\trace}{trace}
\DeclareMathOperator{\rank}{rank}

\begin{document}

\title{Advancing Computational Tools for Analyzing Commutative Hypercomplex Algebras}

\author{
    José Domingo Jiménez López\thanks{Department of Statistics and Operations Research,  University of Ja\'en, Spain} \\
    \texttt{jdomingo@ujaen.es} \\
    \and
    Jes\'us Navarro-Moreno\footnotemark[1] \\
    \texttt{jnavarro@ujaen.es} \\
    \and
    Rosa Mar\'{\i}a Fern\'andez-Alcal\'a\footnotemark[1] \\
    \texttt{rmfernan@ujaen.es}\and
    Juan Carlos Ruiz Molina\footnotemark[1] \\
    \texttt{jcruiz@ujaen.es}
}

\date{}

\maketitle

\begin{abstract}
Commutative hypercomplex algebras offer significant advantages compared to traditional quaternions as they are compatible with linear algebra techniques and efficient computational implementation, which is crucial for their broad applicability. This paper explores a novel family of commutative hypercomplex algebras, designated as $(\alpha\beta)$-tessarines, which extend the system of generalized Segre’s quaternions and, consequently, elliptic quaternions. The core contribution of this paper is the design of theoretical and computational tools for matrices within this algebraic system, including inversion, square root computation, LU factorization with partial pivoting, and determinant calculation. Furthermore, a spectral theory for $(\alpha\beta)$-tessarines is established, covering eigenvalue and eigenvector analysis, the power method, singular value decomposition, rank-$k$ approximation, and the pseudoinverse. Solutions to the classical least squares problem are also provided. These results will  not only advance the fundamental understanding of hypercomplex algebras but also equip researchers with powerful, novel matrix operations that have not been extensively explored in previous related studies.  The theoretical findings are substantiated through real-world examples, including image reconstruction and color face recognition, that highlight the potential of the proposed techniques.

\vskip0.5cm


\end{abstract}

\section{Introduction}\label{Intro}

The study of hypercomplex algebras has expanded significantly in recent years as they can be used to model a wide range of physical and experimental phenomena (see, e.g., \cite{Hahn2016, Valle2021, Vieira2022, ZengSong, Borio, ZhangGao}, and the references therein). Particular attention has been given to the algebras derived from the Cayley-Dickson construction, such as quaternions \cite{Hanson, Took1, Took2, Valle2020, Cariow2021, Guizzo2023, Miron2023, Grassucci2}, dual-quaternions \cite{Grassucci}, octonions \cite{Cariow2023,Bojesomo2024}, sedenions \cite{Bojesomo2024}, and trigintaduonions \cite{Weng2025}.   Although these algebras are non-commutative, they preserve the structure of a normed space. However, as dimensionality increases, they progressively lose essential properties such as associativity and alternativity. This gradual loss of properties explains why most research related to the Cayley-Dickson methodology has primarily focused on quaternion algebra.

Nonetheless, the distinctive properties of tessarines \cite{Pei, Ortolani2017, Cariow2024} have meant that they have attracted significant interest within four-dimensional hypercomplex algebras. Unlike quaternions, tessarines satisfy both commutativity and associativity,   which enhances their compatibility with linear algebra methods and facilitates their computational implementation.

A significant advantage of working with hypercomplex signals   lies in the concept of properness  \cite{Bihan2017,  Nitta2019, Grassucci_proper}. Properness enables hypercomplex signals to be represented in a more compact and efficient manner by reducing the effective dimensionality of the stochastic processes involved. This reduction is achieved by exploiting the specific structural properties in their second-order statistics, leading to substantial computational benefits in various signal processing applications.

In recent years,  parameterized families of hypercomplex algebras have emerged, encompassing a wide range of specific cases depending on particular parameter values. A notable example of these families are generalized Segre's quaternions (GSQ), in which parameter $\alpha$ determines different subfamilies: hyperbolic quaternions  ($\alpha > 0$)  and elliptic quaternions ($\alpha < 0$).  Furthermore, tessarines constitute a specific example of this algebra when $\alpha = -1$. The study presented in \cite{Navarro_Segre} researches GSQ random variables  and signals,  and explores their applications in modeling and forecasting within the framework of GSQ time series. As a particular case of GSQ for $\alpha < 0$,  elliptic quaternion matrices, specifically the spectral properties, have also been extensively studied in \cite{Kosal1}, including fundamental aspects such as the analysis of eigenvalues,  eigenvectors, singular value decomposition (SVD), the pseudoinverse, and solutions to the least squares problem.

Another recently studied parameterized family is that of $\beta$-quaternions, which includes quaternions ($\beta=-1$) and split-quaternions ($\beta=1$) as particular cases.  Within this algebra,   \citep{Navarro_beta_quaternion} introduces  the concepts of first- and second-order properness, 
addresses the problem of minimum mean square error linear estimation of proper $\beta$-quaternion signals,  and proposes adaptive filters under properness conditions.  

One of the main challenges faced by the scientific community when exploring new hypercomplex algebras is the development of toolboxes or computational tools that can enable the implementation of emerging theoretical results. In 2005, Sangwine and Bihan \cite{Sanqwine_LeBihan} developed a MATLAB toolbox to address computational challenges associated with real quaternions and their matrices. This toolbox, which includes image processing capabilities, has undergone continuous improvements, and the latest updated version was released in February 2024. This computational tool has been a key factor in the significant growth of quaternion algebra research when compared to that of other hypercomplex algebras. More recently,  a MATLAB toolbox,  which incorporates image processing methods based on SVD \cite{Kosal2},  has been developed for elliptic quaternion matrices.

The objective of this research is to explore  a new generalized algebra,  named  $(\alpha\beta)$-tessarines,  for $\alpha\in\mathbb{R}-\{0\}$,    $\beta>0$. This algebra  extends GSQ (obtained from $(\alpha\beta)$-tessarines when $\beta=1$) 
and,  consequently,  the elliptic quaternions as studied in \cite{Kosal1} (corresponding to case $\alpha<0$,   $\beta=1$ within the $(\alpha\beta)$-tessarine algebra).    Based on an analysis of the properties of $(\alpha\beta)$-tessarine numbers, this study researches matrix properties that will have significant applications in signal processing, deep learning, machine learning, optimal control systems, image and data compression, and other related fields.   It also establishes a crucial foundation for creating a toolbox that allows the scientific community to utilize the potential of this new algebra.

Specifically, this paper first introduces fundamental concepts and properties for $(\alpha\beta)$-tessarine matrices (Section \ref{Tessarine Domain}). Next, methods for computing the inverse, square root, LU factorization with partial pivoting, and the determinant are presented (Section \ref{Matrix Transformations}). Section \ref{Spectral Theory}  examines spectral theory methods, including eigenvalues and eigenvectors, the power method, SVD, rank-$k$ approximation, and the pseudoinverse. The solution to the least squares problem is discussed in Section \ref{sectiolineal}. The proposed methods are validated with numerical experiments in Section \ref{Examples}. Finally,   Section \ref{section Conclusions} presents concluding remarks, providing   a detailed justification of the $(\alpha\beta)$-tessarine matrix properties proposed,  and discussing the potential applications of these techniques in different areas of statistical signal processing.


\subsection{Notation}\label{section Notation}
This section introduces the basic notation.  Matrices are represented with boldface uppercase letters, column vectors with boldface lowercase letters, and scalar quantities with regular lowercase letters. Superscripts ``$\texttt{T}$'' and ``$*$''  denote the transpose,  and complex transpose, respectively.  The real and imaginary parts of a complex number are denoted by $\Re\{\cdot\}$ and $\Im\{\cdot\}$, respectively.  This notation is also used to represent the real part and  the component corresponding to the imaginary unit $\mathrm{i}$ in an $(\alpha\beta)$-tessarine.

Furthermore, notations $\mathbb{R}$,  $\mathbb{R}^{+}$, $\mathbb{C}$,  $\mathbb{C}_{\alpha}$, $\mathbb{T}_{\alpha\beta}$ and $\mathbb{G}_{\alpha\beta}$ represent the sets of real,  positive real, complex,  $(\alpha)$-complex, $(\alpha\beta)$-tessarines, and generalized $(\alpha\beta)$-tessarines numbers, respectively.  
In this context,   $\mathbf{X}\in \mathbb{L}^{p\times q}$,   for $\mathbb{L} \in \{\mathbb{R}, \mathbb{C}, \mathbb{C}_{\alpha}, \mathbb{T}_{\alpha\beta}, \mathbb{G}_{\alpha\beta}\}$,    denotes a $p\times q$-dimensional matrix with entries in  $\mathbb{L}$,   meaning it is  a real, complex, $(\alpha)$-complex, $(\alpha\beta)$-tessarine, or generalized $(\alpha\beta)$-tessarine matrix.
Similarly $\mathbf{x}\in \mathbb{L}^{p}$  denotes a  $p$-dimensional 
vector  with entries in $\mathbb{L}$.

If $\mathbf{X}\in \mathbb{L}^{p\times p}$  for $\mathbb{L} \in \{\mathbb{R}, \mathbb{C}, \mathbb{C}_{\alpha}, \mathbb{T}_{\alpha\beta}, \mathbb{G}_{\alpha\beta}\}$,  then $\mathbf{X}^{1/2}$ denotes the square root of $\mathbf{X}$, satisfying $\mathbf{X}^{1/2}\mathbf{X}^{1/2}=\mathbf{X}$,   and $\mathbf{X}^{-1}$ denotes the inverse of $\mathbf{X}$,  satisfying $\mathbf{X}^{-1}\mathbf{X}=\mathbf{X}\mathbf{X}^{-1}=\mathbf{I}_{p} $,  where $\mathbf{I}_{p} $  is the identity matrix of dimension $p$. 
The  $(i,j)$th  element of a matrix $\mathbf{X}$ is denoted by $x_{i,j}$. 
Additionally,  $\trace(\cdot)$ represents the trace of a matrix, and  $\mathbf{0}_{p\times q}$ denotes the $p\times q$ zero matrix.  
  
  Finally,  the following notation is used\footnote{This notation  also applies to  $\mathbf{x}\in \mathbb{T}_{\alpha\beta}^p$ or, where  appropriate,  $x \in \mathbb{T}_{\alpha\beta}$.}:

\begin{tabular}{ll}

$\tilde{\mathbf{X}}$ & Associated $(\alpha \beta)$-tessarine of $\mathbf{X}\in \mathbb{T}_{\alpha\beta}^{p\times q}$.\\


$\mathbf{X}^{\texttt{H}_{3-2n}^1}$ &  Hermitian transpose  of  $\mathbf{X}\in \mathbb{T}_{\alpha\beta}^{p\times q}$,    with $n=1$ for $\alpha>0$ and \\ &$n=2$ for $\alpha<0$.\\

$\langle  \cdot ,  \cdot \rangle_n$ & Inner product  of $(\alpha \beta)$-tessarines,  with $n=1$ for $\alpha>0$ and \\ &   $n=2$ for $\alpha<0$.\\ 
$\| \cdot \|_n$ &  Non-homogeneous norm of $(\alpha \beta)$-tessarines,  with $n=1$ for \\ & $\alpha>0$ and   $n=2$ for $\alpha<0$.\\ 


\end{tabular}

\section{$(\alpha\beta)$-Tessarine Domain }\label{Tessarine Domain}

This section establishes  the basic concepts and properties in the $(\alpha\beta)$-tessarine domain.

\begin{defi}
An $(\alpha)$-complex number $x\in \mathbb{C}_{\alpha}$ is defined as $x=a+b\mathrm{i}$, where  $a, b\in \mathbb{R}$,   and  the imaginary unit $\mathrm{i}$ satisfies $\mathrm{i}^2=\alpha$,  with $\alpha\in \mathbb{R}-\{0\}$.
\end{defi}

The $(\alpha)$-complex number system  includes, as particular cases, several important algebras of interest.  Specifically,  values of $\alpha>0$ correspond to hyperbolic complex numbers,  whereas $\alpha<0$ leads to  elliptic complex numbers, and some algebraic properties of these are thoroughly examined in \cite{Kosal1}.   Moreover, for $\alpha=-1$ (i.e.,  $\mathbb{C}_{-1}\equiv\mathbb{C}$ algebra),  the imaginary unit ``$\mathrm{i}$''  is denoted by ``$\epsilon$''.

\begin{defi}\label{defabTessarine}
An $(\alpha\beta)$-tessarine number $x\in \mathbb{T}_{\alpha\beta}$ is defined as $x=a+b\mathrm{i}+c\mathrm{j}+d\mathrm{k}$,  where  $a, b, c, d\in \mathbb{R}$,   and  the imaginary units $(\mathrm{i},\mathrm{j},\mathrm{k})$ satisfy the following multiplication rules:
\begin{equation*}\label{rdf1}
\begin{tabular}{c|ccc|}
    & $\mathrm{i}$ & $\mathrm{j}$ & $\mathrm{k}$ \\\hline
   $\mathrm{i}$ & $\alpha$& $\mathrm{k}$ & $\alpha\mathrm{j}$\\ 
    $\mathrm{j}$ & $\mathrm{k}$ & $\beta$ & $\beta\mathrm{i}$  \\ 
  $ \mathrm{k}$ & $\alpha\mathrm{j}$ & $\beta\mathrm{i}$ &  $\alpha\beta$ \\
  \hline
\end{tabular}
\end{equation*}
where  $\alpha\in \mathbb{R}-\{0\}$ and $\beta \in \mathbb{R}^{+}$.   
\end{defi}

The $(\alpha\beta)$-tessarine numbers contain,  as a particular case,  the algebraic structure of GSQ \cite{Navarro_Segre} when $\beta=1$,   which includes two special subalgebras:  the hyperbolic quaternions  (if $\alpha>0$, $\beta=1$) and the elliptic quaternions  \cite{Kosal1} (if $\alpha<0$, $\beta=1$).   Additionally, when   $\alpha=-1$ and $\beta= 1$,   algebra $\mathbb{T}_{\alpha\beta}$ corresponds to the tessarines   \cite{Vieira2022}. 

According to the multiplication operation,  the $(\alpha\beta)$-tessarine algebra is  commutative, associative, and distributive over addition.  

The algebraic structure of $(\alpha\beta)$-tessarines enables alternative representations to be created, meaning that they can be easily manipulated and analyze.  A key aspect of this framework is the introduction of special units that enable $(\alpha\beta)$-tessarines  to be decomposed into simpler components.  This decomposition not only provides insight into the internal structure of these numbers but also proves useful in various computational and theoretical applications.  The following property formalizes  this decomposition.  

\begin{property}\label{property1}
Consider the following two special units:
\begin{equation}\label{defww}
\mathrm{w}_1=\frac{\sqrt{\beta}+\mathrm{j}}{2\sqrt{\beta}}, \qquad
\mathrm{w}_2=\frac{\sqrt{\beta}-\mathrm{j}}{2\sqrt{\beta}}
\end{equation}
 which satisfy the properties
\begin{equation*}\label{prwww}
\mathrm{w}_1+\mathrm{w}_2=1, \quad \mathrm{w}_1\mathrm{w}_2=0, \qquad \mathrm{w}_{i}^n=\mathrm{w}_{i},\quad  n \in \mathbb{N}\quad i=1,2
\end{equation*}
Then,  the $(\alpha\beta)$-tessarine  $x=a+b\mathrm{i}+c\mathrm{j}+d\mathrm{k}$ can be represented  in the following forms:
\begin{equation}\label{rep}
x=x_a+x_b\mathrm{j}=x_s\mathrm{w}_1+x_d\mathrm{w}_2
\end{equation}
where $x_a=a+b\mathrm{i}$, $x_b=c+d\mathrm{i}$, $x_s=x_a+\sqrt{\beta}x_b$, and  $x_d=x_a-\sqrt{\beta}x_b$ are $(\alpha)$-complex numbers.
\end{property}


For a given $x=a+b\mathrm{i}+c\mathrm{j}+d\mathrm{k}\in  \mathbb{T}_{\alpha\beta}$,  define the function  $s(x)=a+b+c+d$ and introduce the  following auxiliary $(\alpha\beta)$-tessarines: 
\begin{equation}\label{sdf}
\begin{split}
x^{(\theta,\tau)}&=a+\frac{b}{\theta}\mathrm{i}+\frac{c}{\tau}\mathrm{j}+\frac{d}{\tau\theta}\mathrm{k}, \qquad \theta\not =0, \tau>0\\
x^{\mathrm{i}}&=a+b\mathrm{i}-c\mathrm{j}-d\mathrm{k}\\
x^{\mathrm{j}}&=a-b\mathrm{i}+c\mathrm{j}-d\mathrm{k}\\
x^{\mathrm{k}}&=a-b\mathrm{i}-c\mathrm{j}+d\mathrm{k}
\end{split}
\end{equation}
Moreover,  the real vector of $x\in  \mathbb{T}_{\alpha\beta}$ is given  by $\mathbf{x}_r=[a,b,c,d]^{\texttt{T}}$. 



The concept of  the associated  $(\alpha\beta)$-tessarine will play a fundamental role in the subsequent results.  This concept is  formally introduced in the following definition,  for both $\alpha<0$ and $\alpha>0$.

\begin{defi}\label{assov}  Given $x\in  \mathbb{T}_{\alpha\beta}$,  
\begin{description}
\item[ \normalfont{Case }$\alpha<0$:] Let $x=a+c \mathrm{j}$,  with $a,c\in \mathbb{R}$. The  associated $(\alpha\beta)$-tessarine of $x$
is defined as the $(\alpha\beta)$-tessarine $\tilde{x}=\tilde{a}+\tilde{c}\mathrm{j}$, such that  $$x=\frac{s(\tilde{x})}{2}
   +\frac{s(\tilde{x}^{\mathrm{i}})}{2\sqrt{\beta}}\mathrm{j}$$  
    
\item[ \normalfont{Case }$\alpha>0$:]  Let $x=a+b\mathrm{i}+c\mathrm{j}+d\mathrm{k}$,  with $a,b,c,d \in \mathbb{R}$.  The  associated $(\alpha\beta)$-tessarine of $x$ is defined as the $(\alpha\beta)$-tessarine $\tilde{x}=\tilde{a}+\tilde{b}\mathrm{i}+\tilde{c}\mathrm{j}+\tilde{d}\mathrm{k}$, such that  $$x=\frac{s(\tilde{x})}{4}+\frac{s(\tilde{x}^{\mathrm{j}})}{4\sqrt{\alpha}}\mathrm{i}
   +\frac{s(\tilde{x}^{\mathrm{i}})}{4\sqrt{\beta}}\mathrm{j}+\frac{s(\tilde{x}^{\mathrm{k}})}{4\sqrt{\alpha}\sqrt{\beta}}\mathrm{k}$$ 
   \end{description}
   
   Moreover,   $x$ is said to be positive (semipositive), if all the components of its associated  $(\alpha\beta)$-tessarine  $\tilde{x}$  are positive (non-negative).    
   Given $x,y \in \mathbb{T}_{\alpha\beta}$,   the notation $x\succeq y$ or $y\preceq x$ will be used to indicate that  $x-y$ is semipositive.
   \end{defi}
   

\begin{property}\label{plpprop1}
Let  $x_n  \in \mathbb{T}_{\alpha\beta}$,   $n=1,2$.   Then,
\begin{enumerate}
\item $(x_1^{\mathrm{\nu}})^{\mathrm{\nu}}=x_1, \quad \nu=\mathrm{i},\mathrm{j},\mathrm{k}$.
\item $(x_1^{\mathrm{\nu}_1})^{\mathrm{\nu}_2}=x_1^{\mathrm{\nu}_3}, \quad  \nu_1, \nu_2, \nu_3=\mathrm{i},\mathrm{j},\mathrm{k},\quad \nu_1\neq \nu_2 \neq \nu_3 $.
\item $x_1^{\mathrm{\nu}}x_2^{\mathrm{\nu}}=(x_1x_2)^{\mathrm{\nu}}, \quad  \nu=\mathrm{i},\mathrm{j},\mathrm{k}$.
\item $x_1+x_1^{\mathrm{i}}+x_1^{\mathrm{j}}+x_1^{\mathrm{k}}=4 \Re\{x_1\}$.
\item $x_1x_1^{\mathrm{i}}x_1^{\mathrm{j}}x_1^{\mathrm{k}}=\Re\{x_1x_1^{\mathrm{i}}\}^2-\alpha\Im\{x_1x_1^{\mathrm{i}}\}^2$.
\item If $\alpha<0$,  and  $x_n=a_n+c_n\mathrm{j}$ has the associated $(\alpha\beta)$-tessarine $\tilde{x}_n=\tilde{a}_n+\tilde{c}_n\mathrm{j}$,   then  the associated $(\alpha\beta)$-tessarines  of $x_1+x_2$ and  $x_1x_2$
are $(\tilde{a}_1+\tilde{a}_2)+(\tilde{c}_1+\tilde{c}_2)\mathrm{j}$ and $\tilde{a}_1\tilde{a}_2+\tilde{c}_1\tilde{c}_2\mathrm{j}$, respectively.  Moreover, if $\tilde{a}_1,\tilde{c}_1 \in  \mathbb{R}-\{0\}$, then the associated $(\alpha\beta)$-tessarine of $ x_1^{-1}$ is  $\tilde{a}_1^{-1}+\tilde{c}_1^{-1}\mathrm{j}$.
        \item If $\alpha>0$,  and $x_n$
    has the associated $(\alpha\beta)$-tessarine  $\tilde{x}_n=\tilde{a}_n+\tilde{b}_n\mathrm{i}+\tilde{c}_n\mathrm{j}+\tilde{d}_n\mathrm{k}$,   then  the associated $(\alpha\beta)$-tessarines  of $x_1+x_2$ and  $x_1x_2$
are   $(\tilde{a}_1+\tilde{a}_2)+(\tilde{b}_1+\tilde{b}_2)\mathrm{i}+(\tilde{c}_1+\tilde{c}_2)\mathrm{j}+(\tilde{d}_1+\tilde{d}_2)\mathrm{k}$ and $\tilde{a}_1\tilde{a}_2+\tilde{b}_1\tilde{b}_2\mathrm{i}+\tilde{c}_1\tilde{c}_2\mathrm{j}+\tilde{d}_1\tilde{d}_2\mathrm{k}$, respectively.   Moreover, if $\tilde{a}_1,\tilde{b}_1,\tilde{ c}_1, \tilde{d}_1  \in  \mathbb{R}-\{0\}$,  then  the associated $(\alpha\beta)$-tessarine of $ x_1^{-1}$ is $\tilde{a}_1^{-1}+\tilde{b}_1^{-1}\mathrm{i}+\tilde{c}_1^{-1}\mathrm{j}+\tilde{d}_1^{-1}\mathrm{k}$.
\end{enumerate}
\end{property}

All these concepts and properties  can be generalized to the vectorial case.  In particular,
a $p$-dimensional $(\alpha)$-complex  vector $\mathbf{x}\in \mathbb{C}_{\alpha}^{p}$ is given by $\mathbf{x}=[x_{1},\ldots, x_{p}]^\texttt{T}$, 
where $x_{n}=a_{n}+b_{n}\mathrm{i}$,  for $n=1,\ldots, p$,  with $a_{n}, b_{n}\in \mathbb{R}$. 
Similarly,  an $(\alpha)$-complex matrix $\mathbf{X}\in \mathbb{C}_{\alpha}^{p\times q}$ takes the form of  $\mathbf{X}=\mathbf{A}+\mathbf{B}\mathrm{i}$,  with $\mathbf{A}, \mathbf{B}\in \mathbb{R}^{p\times q}$.

Analogously,  a  $p$-dimensional $(\alpha\beta)$-tessarine vector $\mathbf{x}\in \mathbb{T}_{\alpha\beta}^{p}$ is given by $\mathbf{x}=[x_{1},\ldots, x_{p}]^\texttt{T}$,
where $x_{n}=a_{n}+b_{n}\mathrm{i}+c_{n}\mathrm{j}+d_{n}\mathrm{k}$, for $ n=1,\ldots, p$,  with $a_{n}, b_{n}, c_{n}, d_{n}\in \mathbb{R}$.

 The real vector of  $\mathbf{x}\in \mathbb{T}_{\alpha\beta}^{p}$ is denoted as
\begin{equation*}\label{realvect}
\mathbf{x}_r=[\mathbf{a}^{\texttt{T}},\mathbf{b}^{\texttt{T}},\mathbf{c}^{\texttt{T}},\mathbf{d}^{\texttt{T}}]^{\texttt{T}}
\end{equation*} 
where
$\mathbf{a}=[a_1,\ldots,a_p]^\texttt{T}$, $\mathbf{b}=[b_1,\ldots,b_p]^\texttt{T}$, 
$\mathbf{c}=[c_1,\ldots,c_p]^\texttt{T}$, $\mathbf{d}=[d_1,\ldots,d_p]^\texttt{T}$.

From expressions in \eqref{sdf}, the following auxiliaries $(\alpha\beta)$-tessarine vectors are introduced: $\mathbf{x}^{(\theta,\tau)}=[x_1^{(\theta,\tau)}, \ldots,x_p^{(\theta,\tau)}]^\texttt{T}$  and $\mathbf{x}^\nu =[x_1^{\nu}, \ldots,x_p^{\nu}]^\texttt{T}$, for $\nu= {\mathrm{i}},{\mathrm{j}},{\mathrm{k}}$.

Additionally,  an $(\alpha\beta)$-tessarine matrix $\mathbf{X}\in \mathbb{T}_{\alpha\beta}^{p\times q}$ takes the form of  $\mathbf{X}=\mathbf{A}+\mathbf{B}\mathrm{i}+\mathbf{C}\mathrm{j}+\mathbf{D}\mathrm{k}$, with $\mathbf{A}, \mathbf{B},\mathbf{C}, \mathbf{D} \in \mathbb{R}^{p\times q}$.  Analogous to equation \eqref{rep}, $\mathbf{X}\in \mathbb{T}_{\alpha\beta}^{p\times q}$ can be expressed as follows:
\begin{equation}\label{repma}
\mathbf{X}=\mathbf{X}_{a}+\mathbf{X}_{b}\mathrm{j}=\mathbf{X}_{s}\mathrm{w}_1+\mathbf{X}_{d}\mathrm{w}_2
\end{equation}
with $\mathbf{X}_{a}=\mathbf{A}+\mathbf{B}\mathrm{i}$, $\mathbf{X}_{b}=\mathbf{C}+\mathbf{D}\mathrm{i}$, $\mathbf{X}_{s}=\mathbf{X}_a+\sqrt{\beta}\mathbf{X}_b$,  $\mathbf{X}_{d}=\mathbf{X}_a-\sqrt{\beta}\mathbf{X}_b$,  and $\mathrm{w}_n$,  for $n=1,2$,  defined in \eqref{defww}.  Furthermore, in a similar manner to the expressions in \eqref{sdf}, by defining $\mathbf{X}=\left[\mathbf{x}_1,\ldots,\mathbf{x}_q \right]$, with $\mathbf{x}_i\in \mathbb{T}_{\alpha\beta}^{p}$, for $i=1,\ldots,q$, the following auxiliary matrices are introduced:
\begin{equation*}\label{repma2}
\begin{split}
\mathbf{X}^{(\theta,\tau)}&=\mathbf{A}+\frac{\mathbf{B}}{\theta}\mathrm{i}+\frac{\mathbf{C}}{\tau}\mathrm{j}+\frac{\mathbf{D}}{\tau\theta}\mathrm{k}, \qquad \theta\not =0, \tau>0\\
\mathbf{X}^{\nu}&=\left[\mathbf{x}_1^{\nu},\ldots,\mathbf{x}_q^{\nu} \right],\qquad \nu=\mathrm{i},\mathrm{j},\mathrm{k}
\end{split}
\end{equation*}

Extending the concept of the associated $(\alpha\beta)$-tessarine provided in Definition \ref{assov}, matrix $\tilde{\mathbf{X}}$ is defined as the associated $(\alpha\beta)$-tessarine  matrix of $\mathbf{X}$ if  each element $(i,j)$ of $\tilde{\mathbf{X}}$ corresponds to the associated $(\alpha\beta)$-tessarine of element $(i,j)$ of $\mathbf{X}$.  

Furthermore,   the Hermitian transpose with respect to $(\theta, \tau)$ of an $(\alpha\beta)$-tessarine vector  $\mathbf{x}$ (or matrix $\mathbf{X}$)  is denoted as  $\mathbf{x}^{\texttt{H}_\theta^\tau}=\mathbf{x}^{(\theta,\tau)^\texttt{T}}$ (or $\mathbf{X}^{\texttt{H}_{\theta}^\tau}=\mathbf{X}^{(\theta,\tau)^\texttt{T}}$). 

\begin{defi}
 A square matrix $\mathbf{X}\in \mathbb{T}_{\alpha\beta}^{p\times p}$ is said to be $n$-Hermitian,   for $n=1,2$,  if $\mathbf{X}=\mathbf{X}^{\texttt{H}_{3-2n}^{1}}$.
\end{defi}


The following properties detail fundamental algebraic  characteristics of the  $(\alpha\beta)$-tessarine matrices.

 \begin{property}\label{proherm1}
Let  $\mathbf{X}_1\in \mathbb{T}^{p\times q}_{\alpha\beta}$  and  $\mathbf{X}_2\in \mathbb{T}^{q\times r}_{\alpha\beta}$.    Then,  the following properties hold:
\begin{enumerate}
 \item If  $\alpha<0$,  suppose that $\mathbf{X}_n=\mathbf{A}_n+\mathbf{C}_n\mathrm{j}$, for $n=1,2$, have the associated $(\alpha\beta)$-tessarine matrices $\tilde{\mathbf{X}}_n=\tilde{\mathbf{A}}_n+\tilde{\mathbf{C}}_n\mathrm{j}$, respectively. Then,
  the associated    $(\alpha\beta)$-tessarine matrix of $\mathbf{X}_1\mathbf{X}_2$ is given by  $\tilde{\mathbf{A}}_1\tilde{\mathbf{A}}_2+\tilde{\mathbf{C}}_1\tilde{\mathbf{C}}_2\mathrm{j}$. Moreover, if $p=q$, and  the associated $(\alpha\beta)$-tessarine of $\trace(\mathbf{X}_1)$ is denoted by $\tilde{x}$,  it holds that $\tilde{x}=\trace(\tilde{\mathbf{X}}_1)$.
\item If  $\alpha>0$, suppose that $\mathbf{X}_n=\mathbf{A}_n+\mathbf{B}_n\mathrm{i}+\mathbf{C}_n\mathrm{j}+\mathbf{D}_n\mathrm{k}$, for $n=1,2$, have  the associated $(\alpha\beta)$-tessarine matrices  $\tilde{\mathbf{X}}_n=\tilde{\mathbf{A}}_n+\tilde{\mathbf{B}}_n\mathrm{i}+\tilde{\mathbf{C}}_n\mathrm{j}+\tilde{\mathbf{D}}_n\mathrm{k}$, respectively.  Then,
    the associated $(\alpha\beta)$-tessarine matrix of $\mathbf{X}_1\mathbf{X}_2$ is given by  $\tilde{\mathbf{A}}_1\tilde{\mathbf{A}}_2+\tilde{\mathbf{B}}_1\tilde{\mathbf{B}}_2\mathrm{i}+\tilde{\mathbf{C}}_1\tilde{\mathbf{C}}_2\mathrm{j}
    +\tilde{\mathbf{D}}_1\tilde{\mathbf{D}}_2\mathrm{k}$.  Moreover, if $p=q$ and  the associated $(\alpha\beta)$-tessarine of $\trace(\mathbf{X}_1)$ is denoted by $\tilde{x}$, it holds that $\tilde{x}=\trace(\tilde{\mathbf{X}}_1)$.
     \item $(\mathbf{X}_1 \mathbf{X}_2)^{\texttt{H}_{3-2n}^{1}}= \mathbf{X}_2^{\texttt{H}_{3-2n}^{1}}\mathbf{X}_1^{\texttt{H}_{3-2n}^{1}}$,  for $ n=1,2$.
     
   \item  $\mathbf{X}\in \mathbb{T}_{\alpha\beta}^{p\times p}$  is $1$-Hermitian $\Leftrightarrow$ $\mathbf{A}=\mathbf{A}^\texttt{T}$, $\mathbf{B}=\mathbf{B}^\texttt{T}$, $\mathbf{C}=\mathbf{C}^\texttt{T}$,  and $\mathbf{D}=\mathbf{D}^\texttt{T}$.
         \item  $\mathbf{X}\in \mathbb{T}_{\alpha\beta}^{p\times p}$  is $2$-Hermitian $\Leftrightarrow$ $\mathbf{A}=\mathbf{A}^\texttt{T}$, $\mathbf{B}=-\mathbf{B}^\texttt{T}$, $\mathbf{C}=\mathbf{C}^\texttt{T}$,  and $\mathbf{D}=-\mathbf{D}^\texttt{T}$. Moreover, under this condition,  it can be verified that  the  $\mathrm{i}$ and $\mathrm{k}$ parts of $\mathbf{x}^{\texttt{H}_{-1}^{1}}\mathbf{X}\mathbf{x}$ vanish,    $\forall \mathbf{x}\in \mathbb{T}^{p}_{\alpha\beta}$.
      \end{enumerate}
      \end{property}

Several key concepts need to be established within the mathematical framework of $(\alpha\beta)$-tessarine matrices.

\begin{defi}
A square matrix $\mathbf{X}\in \mathbb{T}^{p_\times p}_{\alpha\beta}$ is said to be  positive definite if, and only if,  the following  conditions are satisfied, for $n = 1$ (if  $\alpha > 0$)  and $n = 2$ (if $\alpha < 0$):
\begin{enumerate}
\item $\mathbf{X}$ is $n$-Hermitian.
\item For any nonzero vector $\mathbf{x}\in \mathbb{T}^{p}_{\alpha\beta}$,  $\mathbf{x}^{\texttt{H}_{3-2n}^{1}}\mathbf{X}\mathbf{x}$ is a positive  $(\alpha\beta)$-tessarine number.
\end{enumerate}
\end{defi}

\begin{defi}
The rank of an $(\alpha\beta)$-tessarine matrix $\mathbf{X}\in \mathbb{T}^{p_\times q}_{\alpha\beta}$, denoted by $\rank(\mathbf{X})$,  is defined as the number of independent columns.
\end{defi}

\begin{defi}\label{toeplitz}
A  Toeplitz matrix  is a square matrix $\mathbf{X}\in \mathbb{T}^{p_\times p}_{\alpha\beta}$ whose elements  $x_{i,j}$ satisfy the condition
$x_{i,j}=x_{i+1,j+1}$.
\end{defi}

\begin{defi}
Consider  a square matrix $\mathbf{X}\in \mathbb{T}_{\alpha\beta}^{p\times p}$. The determinant of $\mathbf{X}$, denoted by $\det(\mathbf{X})$,  is defined as
\begin{equation}\label{deterntes}
\det(\mathbf{X})=\sum_{\tau \in P_p}\mathrm{sgn}(\tau)\prod_{i=1}^px_{i,\tau(i)}
\end{equation}
with $P_p$ the set of permutations of $\{1,2,\ldots,p\}$,   $\mathrm{sgn}(\tau)$ the signature of $\tau$,  and $x_{i,\tau(i)}$ is the $(i,\tau(i))$th  entry of $\mathbf{X}$  after the permutation  $\tau$ is applied to the column indices.

\end{defi}

By applying the  commutative, associative and distributive properties of the $(\alpha\beta)$-tessarine product   and Property \ref{plpprop1},  the following result is obtained.
\begin{proposition}\label{porpinv}
Consider  $\mathbf{X}, \mathbf{Y} \in \mathbb{T}_{\alpha\beta}^{p\times p}$. Then,
\begin{enumerate}
\item $\det(\mathbf{X}\mathbf{Y})=\det(\mathbf{X})\det(\mathbf{Y})$.
\item $\det(\mathbf{X}^\nu)=\det(\mathbf{X})^\nu, \quad \nu=\mathrm{i},\mathrm{j},\mathrm{k}$.
\end{enumerate}

Moreover, if $\mathbf{X}$ and $\mathbf{Y}$ are invertible then,
\begin{itemize}
\item[$3_.$] $\det(\mathbf{X}^{-1})=\det(\mathbf{X})^{-1}$.
\item[$4_.$] $(\mathbf{X}^\nu)^{-1}=(\mathbf{X}^{-1})^\nu, \quad \nu=\mathrm{i},\mathrm{j},\mathrm{k}$.
\item[$5_.$] $(\mathbf{X}\mathbf{Y})^{-1}=\mathbf{Y}^{-1}\mathbf{X}^{-1}$.
\end{itemize}
\end{proposition}

Next, given two $(\alpha\beta)$-tessarine matrices $\mathbf{X}, \mathbf{Y} \in \mathbb{T}_{\alpha\beta}^{p\times q}$,  the following product is defined:
\begin{equation}\label{producttoma}
<\mathbf{X},\mathbf{Y}>_{n}=\trace(\mathbf{X}^{\texttt{H}_{3-2n}^{1}}\mathbf{Y})
\end{equation}
with $n=1$ for $\alpha>0$, and  $n=2$ for $\alpha<0$.
\begin{property}\label{rodproxpro}
Let  $\mathbf{X}, \mathbf{Y} \in \mathbb{T}_{\alpha\beta}^{p\times q}$. Then,  the relations below hold:
\begin{enumerate}
\item For $\alpha<0$,   $<\mathbf{X},\mathbf{Y}>_2=<\mathbf{Y},\mathbf{X}>^{\mathrm{j}}_2=<\mathbf{X}^{\texttt{H}_{-1}^1},\mathbf{Y}^{\texttt{H}_{-1}^1}>_2^{\mathrm{j}}$,   with $x^{\mathrm{j}}$ defined in \eqref{sdf}.
\item For $\alpha>0$,   $<\mathbf{X},\mathbf{Y}>_1=<\mathbf{Y},\mathbf{X}>_1=<\mathbf{X}^{\texttt{H}_1^1},\mathbf{Y}^{\texttt{H}_1^1}>_1$.
\end{enumerate}
\end{property}

From the product defined in \eqref{producttoma}, the following non-homogeneous norms are introduced:
\begin{equation}\label{normim}
\|\mathbf{X}\|_n=\Re\{<\mathbf{X},\mathbf{X}>_{n}\}^{1/2}
\end{equation}
with $n=1$ for $\alpha>0$,  and  $n=2$ for $\alpha<0$.
\begin{property}\label{realnrrr}
Let  $x\in \mathbb{T}_{\alpha\beta}$.  Then,
\begin{enumerate}
 \item If $\alpha<0$ and $x=a+c\mathrm{j}$ has the associated $(\alpha\beta)$-tessarine $\tilde{x}=\tilde{a}+\tilde{c}\mathrm{j}$,  then
\begin{equation*}
\|x\|_2^{2}=\frac{\tilde{a}^2+\tilde{c}^2}{2}
\end{equation*}
 \item If $\alpha>0$ and $x=a+b\mathrm{i}+c\mathrm{j}+d\mathrm{k}$ has the associated $(\alpha\beta)$-tessarine $\tilde{x}=\tilde{a}+\tilde{b}\mathrm{i}+\tilde{c}\mathrm{j}+\tilde{d}\mathrm{k}$,   then
\begin{equation*}
\|x\|_1^{2}=\frac{\tilde{a}^2+\tilde{b}^2+\tilde{c}^2+\tilde{d}^2}{4}
\end{equation*}
\end{enumerate}
\end{property}

To conclude,  an extension of the $(\alpha\beta)$-tessarine algebra to an eight-dimensional algebra is introduced, and it is referred to as the generalized $(\alpha\beta)$-tessarine algebra.
\begin{defi}\label{definition generalized tessarine}
A generalized $(\alpha\beta)$-tessarine $x\in \mathbb{G}_{\alpha\beta}$,  with $\alpha,\beta>0$, is defined as 
$$x=a+b\epsilon+c \mathrm{i}+d\epsilon\mathrm{i}+e\mathrm{j}+f\epsilon\mathrm{j}+g \mathrm{k}+h\epsilon\mathrm{k}$$
  where  $a, b, c, d, e,f,g,h\in \mathbb{R}$,   and the imaginary unit $\epsilon$   satisfies  that $\epsilon^2=-1$ and $\epsilon\nu=\nu\epsilon$,  for  $\nu=\mathrm{i},\mathrm{j},\mathrm{k}$.
\end{defi}

Consistent with Property \ref{property1},  the generalized $(\alpha\beta)$-tessarine $x\in \mathbb{G}_{\alpha\beta}$ can be represented in the following forms:
\begin{equation}\label{gentes}
x=x_1+x_2\epsilon=x_3+x_4\mathrm{j}=x_s\mathrm{w}_1+x_d\mathrm{w}_2
\end{equation}
with  $x_1,x_2 \in \mathbb{T}_{\alpha\beta}$ such that
$x_1=a+c \mathrm{i}+e\mathrm{j}+g \mathrm{k}$, $x_2=b+d \mathrm{i}+f\mathrm{j}+h \mathrm{k}$, $x_3=a+b\epsilon+c\mathrm{i}+d\epsilon\mathrm{i}$, $x_4=e+f\epsilon+g\mathrm{i}+h\epsilon\mathrm{i}$, $x_s=x_3+\sqrt{\beta}x_4$, $x_d=x_3-\sqrt{\beta}x_4$,   and $\mathrm{w}_n$, $n=1,2$,  defined in \eqref{defww}. 

\begin{remark}\label{gggff}
The product generated by units $(1,\epsilon,\mathrm{i}, \epsilon\mathrm{i})$ is equivalent to the $(-1\alpha)$-tessarine product.  Moreover,  given the hypernumber $x=a+b\epsilon+c\mathrm{i}+d\epsilon\mathrm{i}$,  its equivalent $(-1\alpha)$-tessarine is defined as
$\breve{x}=a+b\mathrm{i}+c\mathrm{j}+d\mathrm{k}$.
\end{remark}

From Property \ref{plpprop1},  the following properties hold.
%
%
%

\begin{property}
Let  $x=x_1+x_2\epsilon, y=y_1+y_2\epsilon\in \mathbb{G}_{\alpha\beta}$,  with $x_n=a_n+b_n \mathrm{i}+c_n\mathrm{j}+d_n \mathrm{k}$ and  $y_n=e_n+f_n \mathrm{i}+g_n\mathrm{j}+h_n \mathrm{k}$,   for $n=1,2$.  Consider that the  associated $(\alpha\beta)$-tessarines of $x_n$ and  $y_n$  are $\tilde{x}_n=\tilde{a}_n+\tilde{b}_n \mathrm{i}+\tilde{c}_n\mathrm{j}+\tilde{d}_n \mathrm{k}$  and  $\tilde{y}_n=\tilde{e}_n+\tilde{f}_n \mathrm{i}+\tilde{g}_n\mathrm{j}+\tilde{h}_n \mathrm{k}$,  respectively.   Then,
     \begin{enumerate}
     \item The inverse of $x_1$ is given by $$x^{-1}_1=w_1+w_2\epsilon$$   where the associated $(\alpha\beta)$-tessarines of $w_n$,  for $n=1,2$, are $\tilde{w}_1=\frac{\tilde{a}_1}{\tilde{a}_1^2+\tilde{a}_2^2}+\frac{\tilde{b}_1}{\tilde{b}_1^2+\tilde{b}_2^2}\mathrm{i}+\frac{\tilde{c}_1}{\tilde{c}_1^2+\tilde{c}_2^2}\mathrm{j}
         +\frac{\tilde{d}_1}{\tilde{d}_1^2+\tilde{d}_2^2}\mathrm{k}$ and $\tilde{w}_2=-\frac{\tilde{a}_2}{\tilde{a}_1^2+\tilde{a}_2^2}-\frac{\tilde{b}_2}{\tilde{b}_1^2+\tilde{b}_2^2}\mathrm{i}-\frac{\tilde{c}_2}{\tilde{c}_1^2+\tilde{c}_2^2}\mathrm{j}
         -\frac{\tilde{d}_2}{\tilde{d}_1^2+\tilde{d}_2^2}\mathrm{k}$.
    \item The product $xy$ is given by $$xy=l_1+l_2\epsilon$$  where the associated $(\alpha\beta)$-tessarines of $l_n$, for $n=1,2$, are
     $\tilde{l}_1=\Re\{m_1n_1\}+\Re\{m_2n_2\}\mathrm{i}+\Re\{m_3n_3\}\mathrm{j}
         +\Re\{m_4n_4\}\mathrm{k}$ and   $\tilde{l}_2=\Im\{m_1n_1\}+\Im\{m_2n_2\}\mathrm{i}+\Im\{m_3n_3\}\mathrm{j}
         +\Im\{m_4n_4\}\mathrm{k}$,   with
         \begin{equation*}
         \begin{split}
         m_1&=\tilde{a}_1+\tilde{a}_2\epsilon, \  m_2=\tilde{b}_1+\tilde{b}_2\epsilon,  \ m_3=\tilde{c}_1+\tilde{c}_2\epsilon,  \  m_4=\tilde{d}_1+\tilde{d}_2\epsilon\\
          n_1&=\tilde{e}_1+\tilde{e}_2\epsilon,  \  n_2=\tilde{f}_1+\tilde{f}_2\epsilon,  \  n_3=\tilde{g}_1+\tilde{g}_2\epsilon,  \  n_4=\tilde{h}_1+\tilde{h}_2\epsilon
     \end{split}
     \end{equation*}
     \end{enumerate}

\end{property}

The  above definition and properties can be naturally extended to the matrix setting.  Specifically,  
by analogy with Definition \ref{definition generalized tessarine},   
a  generalized $(\alpha\beta)$-tessarine matrix $\mathbf{X}\in \mathbb{G}_{\alpha\beta}^{p\times q}$  takes the form  $$\mathbf{X}=\mathbf{A}+\mathbf{B}\epsilon+\mathbf{C} \mathrm{i}+\mathbf{D}\epsilon\mathrm{i}+\mathbf{E}\mathrm{j}+\mathbf{F}\epsilon\mathrm{j}+\mathbf{G} \mathrm{k}+\mathbf{H}\epsilon\mathrm{k}$$  with $\mathbf{A}, \mathbf{B}, \mathbf{C}, \mathbf{D}, \mathbf{E}, \mathbf{F}, \mathbf{G},\mathbf{H} \in \mathbb{R}^{p\times q}$.

Furthermore, similarly to equation \eqref{gentes},   $\mathbf{X}\in \mathbb{G}_{\alpha\beta}^{p\times q}$ 
can be alternatively expressed as follows:
\begin{equation}\label{repmaG}
\mathbf{X}=\mathbf{X}_{1}+\mathbf{X}_{2}\epsilon=\mathbf{X}_{3}+\mathbf{X}_{4}\mathrm{j}=\mathbf{X}_{s}\mathrm{w}_1+\mathbf{X}_{d}\mathrm{w}_2
\end{equation}
where  $\mathbf{X}_1, \mathbf{X}_2\in \mathbb{T}_{\alpha\beta}^{p\times q}$, with $\mathbf{X}_1=\mathbf{A}+\mathbf{C}\mathrm{i}+\mathbf{E} \mathrm{j}+\mathbf{G}\mathrm{k}$ and  $\mathbf{X}_2=\mathbf{B}+\mathbf{D}\mathrm{i}+\mathbf{F} \mathrm{j}+\mathbf{H}\mathrm{k}$,  $\mathbf{X}_{3}=\mathbf{A}+\mathbf{B}\epsilon+\mathbf{C} \mathrm{i}+\mathbf{D}\epsilon\mathrm{i}$, $\mathbf{X}_{4}=\mathbf{E}+\mathbf{F}\epsilon+\mathbf{G} \mathrm{i}+\mathbf{H}\epsilon\mathrm{i}$, $\mathbf{X}_{s}=\mathbf{X}_3+\sqrt{\beta}\mathbf{X}_4$,  $\mathbf{X}_{d}=\mathbf{X}_3-\sqrt{\beta}\mathbf{X}_4$,  and $\mathrm{w}_n$, $n=1,2$,  defined in \eqref{defww}.

\section{Methods for $(\alpha\beta)$-Tessarine Matrix Transformations}\label{Matrix Transformations}

This section discusses various computational techniques for performing fundamental transformations on $(\alpha\beta)$-tessarine matrices. Specifically,  they include matrix inversion, matrix square root,   LU factorization with partial pivoting, and determinant computation.  The methodology used to derive the results in this section aligns with that presented in \cite{Kosal1}.  This methodology relies on the results for $(\alpha)$-complex matrices,  which are described in  \ref{Ap1}.

\subsection{Inverse  of an  $(\alpha\beta)$-Tessarine Matrix}
This section presents an efficient approach for computing the inverse of an $(\alpha\beta)$-tessarine matrix. 

\begin{proposition}\label{propoiinvt}
Consider $\mathbf{X}\in \mathbb{T}_{\alpha\beta}^{p\times p}$ expressed in the form \eqref{repma} as $\mathbf{X}=\mathbf{X}_{s}\mathrm{w}_1+\mathbf{X}_{d}\mathrm{w}_2$.  Then,  $\mathbf{X}$ is invertible if, and only if,  matrices $\mathbf{X}_{s}$  and $\mathbf{X}_{d}$ are both invertible.  


Moreover, if $\mathbf{X}$ is invertible, its inverse is given by 
$\mathbf{X}^{-1}=\mathbf{C}+\mathbf{D}\mathrm{j}$,
with $\mathbf{C}=\frac{\mathbf{X}_s^{-1}+\mathbf{X}_d^{-1}}{2}$ and  $\mathbf{D}=\frac{\mathbf{X}_s^{-1}-\mathbf{C}}{\sqrt{\beta}}$,  where
$\mathbf{X}_s^{-1}$ and $\mathbf{X}_d^{-1}$ are obtained from Lemma \ref{lema1}.
\end{proposition}

As an extension, a method for computing the inverse of a generalized  $(\alpha\beta)$-tessarine matrix is proposed below.
\begin{proposition}\label{invtr}
Consider $\mathbf{X}\in \mathbb{G}_{\alpha\beta}^{p\times p}$  expressed in the form  \eqref{repmaG} as $\mathbf{X}=\mathbf{X}_{s}\mathrm{w}_1+\mathbf{X}_{d}\mathrm{w}_2$. Then,   $\mathbf{X}^{-1}$ is invertible  if, and only if, $\mathbf{X}_{s}$  and $\mathbf{X}_{d}$ are both invertible.  


Moreover,    if $\mathbf{X}$ is invertible, its inverse is given by  $\mathbf{X}^{-1}=\mathbf{C}+\mathbf{D}\mathrm{j}$
with $\mathbf{C}=\frac{\mathbf{X}_s^{-1}+\mathbf{X}_d^{-1}}{2}$ and  $\mathbf{D}=\frac{\mathbf{X}_s^{-1}-\mathbf{C}}{\sqrt{\beta}}$.
\end{proposition}
\begin{remark}
Given the equivalence between the $(-1\alpha)$-tessarine product and the one generated by $(1, \epsilon, \mathrm{i}, \epsilon \mathrm{i})$ (see Remark \ref{gggff}),  $\mathbf{X}_s^{-1}$ and $\mathbf{X}_d^{-1}$ can be  obtained by calculating their equivalent $(-1\alpha)$-tessarine matrices, $\breve{\mathbf{X}}_s^{-1}$ and $\breve{\mathbf{X}}_d^{-1}$,  using Proposition \ref{propoiinvt}.
\end{remark}

\begin{corollary}\label{injnj}$\quad$
\begin{enumerate}
\item Consider $\mathbf{X}=\frac{\mathbf{A}_1+\mathbf{A}_2}{2}+\frac{\mathbf{A}_1-\mathbf{A}_2}{2\sqrt{\beta}}\mathrm{j} \in \mathbb{T}_{\alpha\beta}^{p\times p} $,  where $\mathbf{A}_n=\mathbf{M}_n+\mathbf{N}_n\mathrm{i}$, $n=1,2$. Suppose that $\mathbf{X}^{-1}$ exists. Then,
    \begin{equation*}
    \mathbf{X}^{-1}=\frac{\mathbf{A}_1^{-1}+\mathbf{A}_2^{-1}}{2}+\frac{\mathbf{A}_1^{-1}-\mathbf{A}_2^{-1}}{2\sqrt{\beta}}\mathrm{j}
    \end{equation*}
    \item Consider
    \begin{equation*}
    \begin{split}
    \mathbf{X}=&\frac{\mathbf{A}_1+\mathbf{A}_2+\mathbf{A}_3+\mathbf{A}_4}{4}+\frac{\mathbf{A}_1-\mathbf{A}_2+\mathbf{A}_3-\mathbf{A}_4}{4\sqrt{\alpha}}\mathrm{i}\\
    &+\left[\frac{\mathbf{A}_1+\mathbf{A}_2-\mathbf{A}_3-\mathbf{A}_4}{4\sqrt{\beta}}+\frac{\mathbf{A}_1-\mathbf{A}_2-\mathbf{A}_3+\mathbf{A}_4}{4\sqrt{\alpha}\sqrt{\beta}}\mathrm{i}\right]\mathrm{j} \in \mathbb{G}_{\alpha\beta}^{p\times p}
    \end{split}
    \end{equation*}
   where $\mathbf{A}_n=\mathbf{M}_n+\mathbf{N}_n\epsilon$, $n=1,\ldots,4$.  Suppose that $\mathbf{X}^{-1}$ exists. Then,
        \begin{equation*}
    \begin{split}
    \mathbf{X}^{-1}=&\frac{\mathbf{A}_1^{-1}+\mathbf{A}_2^{-1}+\mathbf{A}_3^{-1}+\mathbf{A}_4^{-1}}{4}+\frac{\mathbf{A}_1^{-1}-\mathbf{A}_2^{-1}+\mathbf{A}_3^{-1}-\mathbf{A}_4^{-1}}{4\sqrt{\alpha}}\mathrm{i}\\
    &+\left[\frac{\mathbf{A}_1^{-1}+\mathbf{A}_2^{-1}-\mathbf{A}_3^{-1}-\mathbf{A}_4^{-1}}{4\sqrt{\beta}}+\frac{\mathbf{A}_1^{-1}-\mathbf{A}_2^{-1}-\mathbf{A}_3^{-1}+\mathbf{A}_4^{-1}}{4\sqrt{\alpha}\sqrt{\beta}}\mathrm{i}\right]\mathrm{j}
    \end{split}
    \end{equation*}
\end{enumerate}
\end{corollary}

\subsection{Square  Root of an $(\alpha\beta)$-Tessarine Matrix}
Next,   an efficient method for computing the square root of an  $(\alpha\beta)$-tessarine matrix is provided.
\begin{proposition}\label{proinmb}
Consider $\mathbf{X}\in \mathbb{T}_{\alpha\beta}^{p\times p}$  expressed in form \eqref{repma} as $\mathbf{X}=\mathbf{X}_{s}\mathrm{w}_1+\mathbf{X}_{d}\mathrm{w}_2$.  Then,  $\mathbf{X}$ possesses a square root $\mathbf{X}^{1/2}$ if, and only if,  the square roots of both $\mathbf{X}_{s}$ and $\mathbf{X}_{d}$,  $\mathbf{X}_{s}^{1/2}$ and $\mathbf{X}_{d}^{1/2}$,   exist.


Moreover,   if $\mathbf{X}^{1/2}$ exists,  it is given by
$\mathbf{X}^{1/2}=\mathbf{C}+\mathbf{D}\mathrm{j}$,
with $\mathbf{C}=\frac{\mathbf{X}_s^{1/2}+\mathbf{X}_d^{1/2}}{2}$ and  $\mathbf{D}=\frac{\mathbf{X}_s^{1/2}-\mathbf{C}}{\sqrt{\beta}}$, where
$\mathbf{X}_s^{1/2}$ and $\mathbf{X}_d^{1/2}$ are obtained from Lemma \ref{lema221}.

\end{proposition}
\begin{remark}
Matrix $\mathbf{X}^{1/2}$ is a generalized  $(\alpha\beta)$-tessarine matrix  and may  not be unique.
\end{remark}

To conclude this section,  a method for computing the square root of a generalized   $(\alpha\beta)$-tessarine matrix is provided.
\begin{proposition}
Consider $\mathbf{X}\in \mathbb{G}_{\alpha\beta}^{p\times p}$  expressed in form  \eqref{repmaG} as $\mathbf{X}=\mathbf{X}_{s}\mathrm{w}_1+\mathbf{X}_{d}\mathrm{w}_2$. Then,   $\mathbf{X}$ possesses a square root $\mathbf{X}^{1/2}$ if, and only if,  the square roots of both $\mathbf{X}_{s}$ and $\mathbf{X}_{d}$,  $\mathbf{X}_{s}^{1/2}$ and $\mathbf{X}_{d}^{1/2}$,   exist.

Moreover,   if $\mathbf{X}^{1/2}$ exists,  it is given by
$\mathbf{X}^{1/2}=\mathbf{C}+\mathbf{D}\mathrm{j}$,
with $\mathbf{C}=\frac{\mathbf{X}_s^{1/2}+\mathbf{X}_d^{1/2}}{2}$ and  $\mathbf{D}=\frac{\mathbf{X}_s^{1/2}-\mathbf{C}}{\sqrt{\beta}}$.  
\end{proposition}

\begin{remark}
Given the equivalence between the  $(-1\alpha)$-tessarine product and the one generated by  $(1,\epsilon,\mathrm{i}, \epsilon\mathrm{i})$ (see Remark \ref{gggff}), square roots 
$\mathbf{X}_s^{1/2}$ and $\mathbf{X}_d^{1/2}$ are obtained  by computing their corresponding  $(-1\alpha)$-tessarine matrices, $\breve{\mathbf{X}}_s^{1/2}$ and $\breve{\mathbf{X}}_d^{1/2}$,  from Proposition \ref{proinmb}.  
\end{remark}

\subsection{LU Factorization with Partial Pivoting of an $(\alpha\beta)$-Tessarine Matrix}
 The problem of LU factorization with partial pivoting of a complex matrix consists of determining the matrices $(\mathbf{P}, \mathbf{L}, \mathbf{U})$ that satisfy decomposition  $\mathbf{P}\mathbf{X}=\mathbf{L}\mathbf{U}$,  where $\mathbf{P}$
 is a permutation matrix,  and $\mathbf{L}$ and $\mathbf{U}$ are the lower and upper triangular matrices, respectively.   Next,  this methodology is extended  to the  $(\alpha\beta)$-tessarine domain.

\begin{proposition}\label{LUprpp}
Consider $\mathbf{X}\in \mathbb{T}_{\alpha\beta}^{p\times p}$ is expressed in form  \eqref{repma} as    $\mathbf{X}=\mathbf{X}_{s}\mathrm{w}_1+\mathbf{X}_{d}\mathrm{w}_2$.   Let $(\mathbf{P}_s,\mathbf{L}_s,\mathbf{U}_s)$ and  $(\mathbf{P}_d,\mathbf{L}_d,\mathbf{U}_d)$ be the matrices of the LU factorizations with partial pivoting  of the $(\alpha)$-complex matrices $\mathbf{X}_{s}$ and $\mathbf{X}_{d}$,  respectively,  obtained via Lemma \ref{LUlema}. Then, the matrices $(\mathbf{P},\mathbf{L},\mathbf{U})$ of the LU factorization with partial pivoting  of matrix $\mathbf{X}$ take the form
        \begin{equation}\label{matriceslluupp}
        \begin{split}
        \mathbf{P}&=\frac{\mathbf{P}_s+\mathbf{P}_d}{2}+\frac{\mathbf{P}_s-\mathbf{P}_d}{2\sqrt{\beta}}\mathrm{j}\\
         \mathbf{L}&=\frac{\mathbf{L}_s+\mathbf{L}_d}{2}+\frac{\mathbf{L}_s-\mathbf{L}_d}{2\sqrt{\beta}}\mathrm{j}\\
          \mathbf{U}&=\frac{\mathbf{U}_s+\mathbf{U}_d}{2}+\frac{\mathbf{U}_s-\mathbf{U}_d}{2\sqrt{\beta}}\mathrm{j}
        \end{split}
        \end{equation}
\end{proposition}


\begin{remark}
Unlike the complex case,  matrix $\mathbf{P}$  is not necessarily  a binary permutation matrix.   However, it is ortogonal,  i.e., 
$\mathbf{P}^{\texttt{H}_{3-2n}^{1}}\mathbf{P}=\mathbf{P}\mathbf{P}^{\texttt{H}_{3-2n}^{1}}=\mathbf{I}_p$,  where $n=1$ for $\alpha>0$,
and $n=2$ for $\alpha<0$.
\end{remark}

A noteworthy consequence of Proposition \ref{LUprpp} is the efficient computation of the determinant of a matrix, as described in equation  \eqref{deterntes}.  In particular,  the following Corollary provides a formula for computing the determinant of an $(\alpha\beta)$-tessarine matrix,  based  on  its LU factorization with partial pivoting.  It expresses the determinant in terms of the determinant of matrix $\mathbf{P}$ and the diagonal elements of the upper triangular matrix $\mathbf{U}$ in the factorization.  For the computation of the determinant of $\mathbf{P}$,  the following matrices and vectors have been considered: 
\begin{itemize}
 \item  If $\alpha<0$,  let $\mathbf{P}_s$ and  $\mathbf{P}_d$ denote the real  matrices  in Lemma \ref{LUlema} 
 associated with  $\mathbf{X}_{s}$ and $\mathbf{X}_{d}$ given in \eqref{repma},   respectively.   Define vector  $\boldsymbol{\gamma}_{2}=[\gamma_s,\gamma_d]$,  where  $\gamma_s$,  $\gamma_d$  represents the signatures associated with  $\mathbf{P}_s$ and  $\mathbf{P}_d$, respectively. 
 \item  If $\alpha>0$,  let $(\mathbf{P}_{s1}, \mathbf{P}_{d1})$  and  $(\mathbf{P}_{s2}, \mathbf{P}_{d2})$ denote the real  matrices in  \eqref{matrds} associated with $\mathbf{X}_{s}$ and $\mathbf{X}_{d}$ given in  \eqref{repma},  respectively.
Define vector $\boldsymbol{\gamma}_{1}=[\gamma_{s1},\gamma_{d1},\gamma_{s2},\gamma_{d2}]$,  where  $\gamma_{s1},\gamma_{d1},\gamma_{s2},\gamma_{d2}$ 
are the signatures associated with the above four matrices,  respectively.  
 \end{itemize}
Then,  from Propositions \ref{porpinv} and \ref{LUprpp},   the following result holds.

\begin{corollary} Consider $\mathbf{X}\in \mathbb{T}_{\alpha\beta}^{p\times p}$,  and $(\mathbf{P}, \mathbf{L}, \mathbf{U})$ the matrices of the LU factorization with partial pivoting of $\mathbf{X}$ of the form  \eqref{matriceslluupp}. Then,
\begin{enumerate}
\item The determinant of matrix $\mathbf{P}$,  $\det(\mathbf{P})$,   is provided in Table \ref{tabdt} as a  function of the values of $\boldsymbol{\gamma}_1$ ($\alpha>0$) and $\boldsymbol{\gamma}_2$ ($\alpha<0$).
\item The determinant of $\mathbf{X}$ is given by 
\begin{equation*}\label{eficdete}
\det(\mathbf{X})=\det(\mathbf{P})  \prod_{i=1}^p u_{i,i}
\end{equation*}
where $u_{i,i}$  is  the $(i,i)$th element of matrix $\mathbf{U}$. 
\end{enumerate}
\end{corollary}

\begin{remark}
From a computational standpoint,  the calculation of the values of $\boldsymbol{\gamma}_i$,  for $i=1,2$,  
is both efficient and straightforward.
\end{remark}

\begin{table}[h]
\renewcommand{\arraystretch}{1.5} 
\begin{center}
\begin{tabular}{|c|c|}
  \hline
\multicolumn{2}{|c|}{$\alpha<0$} \\ \hline
 $\boldsymbol{\gamma}_{2}$ & $\det(\mathbf{P})$  \\\hline
  $[1,1]$& $1$  \\
  $[-1,-1]$ & $-1$ \\
   $[1,-1]$& $\frac{1}{\sqrt{\beta}}\mathrm{j}$  \\
   $[-1,1]$ & $\frac{-1}{\sqrt{\beta}}\mathrm{j}$  \\
  \hline
\end{tabular}
\hspace*{0.5cm}
\begin{tabular}{|c|c|c|c|}
  \hline
\multicolumn{4}{|c|}{$\alpha>0$} \\ \hline
  $\boldsymbol{\gamma}_{1}$ & $\det(\mathbf{P})$ & 
   $\boldsymbol{\gamma}_{1}$ & $\det(\mathbf{P})$\\\hline
  $[1,1,1,1]$& $1$ & $[-1,-1,-1,-1]$ & $-1$ \\
  $[1,-1,1,-1]$ & $\frac{1}{\sqrt{\alpha}}\mathrm{i}$ & $[-1,1,-1,1]$ & $\frac{-1}{\sqrt{\alpha}}\mathrm{i}$\\
   $[1,1,-1,-1]$  & $\frac{1}{\sqrt{\beta}}\mathrm{j}$ &  $[-1,-1,1,1]$  & $\frac{-1}{\sqrt{\beta}}\mathrm{j}$ \\
 $[1,-1,-1,1]$ &  $\frac{1}{\sqrt{\alpha\beta}}\mathrm{k}$ &  $[-1,1,1,-1]$ &  $\frac{-1}{\sqrt{\alpha\beta}}\mathrm{k}$ \\
   $[1,1,1,-1]$ & $a_1$ &   $[1,-1,-1,-1]$ & $a_2$ \\
   $[1,-1,1,1]$ & $a_1^{\mathrm{i}}$ & $[-1,-1,1,-1]$  &  $a_2^{\mathrm{i}}$ \\
  $[1,1,-1,1]$  & $a_1^{\mathrm{j}}$  &$[-1,1,-1,-1]$  &  $a_2^{\mathrm{j}}$ \\
  $[-1,1,1,1]$ & $a_1^{\mathrm{k}}$ &$[-1,-1,-1,1]$  &  $a_2^{\mathrm{k}}$ \\
  \hline
  \multicolumn{4}{|c|}{$a_1=0.5(1+\frac{1}{\sqrt{\alpha}}\mathrm{i}+\frac{1}{\sqrt{\beta}}\mathrm{j}-\frac{1}{\sqrt{\alpha\beta}}\mathrm{k})$} \\ 
  \multicolumn{4}{|c|}{$a_2=0.5(-1+\frac{1}{\sqrt{\alpha}}\mathrm{i}+\frac{1}{\sqrt{\beta}}\mathrm{j}+\frac{1}{\sqrt{\alpha\beta}}\mathrm{k})$}\\
  \hline
\end{tabular}
\end{center}
\caption{ Determinant of matrix $\mathbf{P}$ as a  function of  $\boldsymbol{\gamma}_{1}$ $(\alpha>0)$ and $\boldsymbol{\gamma}_2$ $(\alpha<0)$.
\label{tabdt}}
\end{table}

%

\section{$(\alpha\beta)$-Tessarine Spectral Theory}\label{Spectral Theory}

This section develops the spectral theory of $(\alpha\beta)$-tessarine matrices, extending fundamental decomposition methods to this algebraic structure.   First,  the eigendecomposition of a square $(\alpha\beta)$-tessarine matrix is established, offering insights into its spectral characteristics.  Next,  the SVD of a general $(\alpha\beta)$-tessarine matrix is derived.   Finally, the pseudoinverse of  an  $(\alpha\beta)$-tessarine matrix is examined, and the SVD approach is employed to handle non-invertible cases. These results are derived from  those established for  $(\alpha)$-complex matrices, which are described in  \ref{Ap1}.

\subsection{Eigendecomposition of an $(\alpha\beta)$-Tessarine Matrix}

The eigendecomposition of an $(\alpha\beta)$-tessarine matrix is now addressed. 
To maintain conciseness,   and where no ambiguity arises,  any eigenvalue and its corresponding eigenvector from eigenvalue set $\{\lambda_i\}_{i=1}^p$ and eigenvector set $\{\mathbf{u}_i\}_{i=1}^p$  are denoted by $\lambda$ and $\mathbf{u}$, respectively.

Following a similar approach to that of \cite{Pei},  a general form of the eigendecomposition for an $(\alpha\beta)$-tessarine matrix is established in Proposition \ref{prpp2}.

\begin{proposition}\label{prpp2}
 Consider $\mathbf{X}\in \mathbb{T}_{\alpha\beta}^{p\times p}$ expressed in form  \eqref{repma} as  $\mathbf{X}=\mathbf{X}_{s}\mathrm{w}_1+\mathbf{X}_{d}\mathrm{w}_2$.  Let $\lambda_{s}$ and $\lambda_{d}$  (respectively, $\mathbf{u}_{s}$ and $\mathbf{u}_{d}$) be the eigenvalues (respectively,  eigenvectors) of $\mathbf{X}_s$ and $\mathbf{X}_d$ computed from Lemma \ref{lema2}. Then, the eigenvalues and eigenvectors of $\mathbf{X}$ are given by $\lambda=\lambda_1+\lambda_2\mathrm{j}$ and $\mathbf{u}=\mathbf{u}_1+\mathbf{u}_2\mathrm{j}$, respectively,  where $\lambda_1=\frac{\lambda_{s}+\lambda_{d}}{2}$, $\lambda_2=\frac{\lambda_{s}-\lambda_1}{\sqrt{\beta}}$, $\mathbf{u}_1=\frac{\mathbf{u}_{s}+\mathbf{u}_{d}}{2}$,   and  $\mathbf{u}_2=\frac{\mathbf{u}_{s}-\mathbf{u}_1}{\sqrt{\beta}}$.
\end{proposition}
\begin{remark}\label{foreig}
If $\alpha<0$,  matrix $\mathbf{X}$ has $p^2$ eigenvalues,  which are  $(\alpha\beta)$-tessarine numbers taking the form
\begin{equation}\label{ealphane}
\lambda=\frac{\Re\{\lambda_s+\lambda_d\}}{2}+\frac{\Im\{\lambda_s+\lambda_d\}}{2}\mathrm{i}
+\frac{\Re\{\lambda_s-\lambda_d\}}{2\sqrt{\beta}}\mathrm{j}+\frac{\Im\{\lambda_s-\lambda_d\}}{2\sqrt{\beta}}\mathrm{k}
\end{equation}

 If $\alpha>0$,  matrix $\mathbf{X}$ has $p^4$ eigenvalues,  which are  generalized $(\alpha\beta)$-tessarine numbers.  Specifically,  let $\{\theta_{ni}\}_{i=1}^p$,  for $n=1,\ldots,4$,  denote the set of eigenvalues corresponding to the following matrices:
 \begin{equation}\label{matricessAA}
 \begin{split}
 \mathbf{A}_1&=\Re\{\mathbf{X}_s\}+\sqrt{\alpha}\Im\{\mathbf{X}_s\}  \qquad  \mathbf{A}_2=\Re\{\mathbf{X}_s\}-\sqrt{\alpha}\Im\{\mathbf{X}_s\} \\
 \mathbf{A}_3&=\Re\{\mathbf{X}_d\}+\sqrt{\alpha}\Im\{\mathbf{X}_d\}  \qquad 
\mathbf{A}_4=\Re\{\mathbf{X}_d\}-\sqrt{\alpha}\Im\{\mathbf{X}_d\}
\end{split}
\end{equation}

For notational simplicity,  $\theta_n$ represents any eigenvalue belonging to set $\{\theta_{ni}\}_{i=1}^p$,  for $n=1,\ldots,4$. Then, the $p^4$ eigenvalues of $\mathbf{X}$ take form
 \eqref{gentes} as:
     \begin{equation}\label{eigenvnv}
     \lambda=\varphi+\psi\epsilon \in \mathbb{G}_{\alpha\beta}
     \end{equation}
 where  the  corresponding $(\alpha\beta)$-tessarines $ \varphi$ and $\psi $ are
  given by   $\varphi= \Re\{\theta_1\}+ \Re\{\theta_2\}\mathrm{i}+ \Re\{\theta_3\}\mathrm{j}+ \Re\{\theta_4\} \mathrm{k} $   and $\psi= \Im\{\theta_1\}+ \Im\{\theta_2\}\mathrm{i}+ \Im\{\theta_3\}\mathrm{j}+ \Im\{\theta_4\} \mathrm{k} $, respectively.

Similar considerations apply to the eigenvectors.
\end{remark}

Next,  an efficient eigenrepresentation of an  $(\alpha\beta)$-tessarine matrix is provided.
 \begin{theorem}\label{prooo3}
  Consider $\mathbf{X}\in \mathbb{T}_{\alpha\beta}^{p\times p}$ expressed in form  \eqref{repma} as  $\mathbf{X}=\mathbf{X}_{s}\mathrm{w}_1+\mathbf{X}_{d}\mathrm{w}_2$.  The following results hold:
\begin{description}
\item[ \normalfont{Case }$\alpha<0$:] Let  $\{\tau_1(\lambda_{si})\}_{i=1}^p$ and $\{\tau_2(\lambda_{di})\}_{i=1}^p$ represent any  two permutations of the sets of eigenvalues of  $\mathbf{X}_s$ and $\mathbf{X}_d$, respectively. Define diagonal matrix $\mathbf{\Lambda}$, such that its $(i,i)$th element $\lambda_i$ is given by the expression \eqref{ealphane}, from $\tau_1(\lambda_{si})$ and $\tau_2(\lambda_{di})$.  

Denote by $\mathbf{W}_n$,  for $n=1,2$,  the matrices whose $p$  columns correspond to the eigenvectors associated with $\{\tau_1(\lambda_{si})\}_{i=1}^p$ and $\{\tau_2(\lambda_{di})\}_{i=1}^p$,  respectively.  Similarly, let $\mathbf{U}$ be the matrix whose $p$ columns correspond to the eigenvectors associated with $\{\lambda_i\}_{i=1}^p$.  Therefore, 
\begin{equation*}
\exists \mathbf{U}^{-1} \Leftrightarrow\exists  \mathbf{W}_n^{-1}, \forall n=1,2
\end{equation*}

Furthermore,  if $\mathbf{U}$ is invertible,  then $\mathbf{X}$ admits the eigenvalue decomposition 
\begin{equation*}\label{eizneg}
\mathbf{X}=\mathbf{U}\mathbf{\Lambda}\mathbf{U}^{-1}
\end{equation*}

 \item[ \normalfont{Case }$\alpha>0$:] Let  $\{\tau_n(\theta_{ni})\}_{i=1}^p$, for $n=1,\ldots,4$,  represent any permutation of the sets of eigenvalues $\{\theta_{ni}\}_{i=1}^p$,  for $n=1,\ldots,4$,  of matrices $\mathbf{A}_n$ defined  in  \eqref{matricessAA}.      Define diagonal matrix $\mathbf{\Lambda}$, such that its $(i,i)$th element $\lambda_i$ is given by the expression \eqref{eigenvnv} from $\tau_n(\theta_{ni})$,  for $n=1,\ldots,4$. 
 
 Denote by $\mathbf{V}_n$,  for $n=1,\dots,4$,  the matrices whose $p$ columns correspond to the eigenvectors associated with $\{\tau_n(\theta_{ni})\}_{i=1}^p$.  Similarly,  let $\mathbf{U}$ be the matrix whose $p$  columns correspond to the eigenvectors associated to $\{\lambda_i\}_{i=1}^p$.  Therefore,  
\begin{equation}\label{hjyuikj}
\exists \mathbf{U}^{-1} \Leftrightarrow\exists  \mathbf{V}_n^{-1}, \forall n=1,\ldots,4
\end{equation}

Furthermore,  if $\mathbf{U}$ is invertible,  then $\mathbf{X}$ admits the eigenvalue decomposition 
\begin{equation}\label{eigdes}
\mathbf{X}=\mathbf{U}\mathbf{\Lambda}\mathbf{U}^{-1}
\end{equation}
\end{description}
 \end{theorem}
 
\begin{proof}
The proof for  case $\alpha>0$ is outlined.  The proof for case $\alpha<0$ follows analogously,  and it has been omitted for brevity.

From Lemma \ref{lema2} and Proposition \ref{prpp2}, it follows that 
    \begin{equation}\label{matV}
    \begin{split}
    \mathbf{U}=&\frac{\mathbf{V}_1+\mathbf{V}_2+\mathbf{V}_3+\mathbf{V}_4}{4}+\frac{\mathbf{V}_1-\mathbf{V}_2+\mathbf{V}_3-\mathbf{V}_4}{4\sqrt{\alpha}}\mathrm{i}\\
    &+\left[\frac{\mathbf{V}_1+\mathbf{V}_2-\mathbf{V}_3-\mathbf{V}_4}{4\sqrt{\beta}}+\frac{\mathbf{V}_1-\mathbf{V}_2-\mathbf{V}_3+\mathbf{V}_4}{4\sqrt{\alpha}\sqrt{\beta}}\mathrm{i}\right]\mathrm{j}
    \end{split}
    \end{equation}

Thus, by applying Lemma \ref{lema1}, Propositions \ref{propoiinvt}  and \ref{invtr}, and Corollary  \ref{injnj},   the result in \eqref{hjyuikj} is established.  Finally,  equation \eqref{eigdes} follows directly  from identity  $\mathbf{X}\mathbf{U}=\mathbf{U}\mathbf{\Lambda}$.
 \end{proof}

\begin{remark}
It is important to highlight that the preceding result offers a simplified and efficient approach for obtaining the spectral representation of $\mathbf{X}$ without requiring the full set of eigenvalues.
\end{remark}

\begin{corollary}\label{Cor3}
 Consider $\mathbf{X}\in \mathbb{T}_{\alpha\beta}^{p\times p}$ expressed in form  \eqref{repma} as  $\mathbf{X}=\mathbf{X}_{s}\mathrm{w}_1+\mathbf{X}_{d}\mathrm{w}_2$.  The following properties hold:
\begin{description}
\item[ \normalfont{Case }$\alpha<0$:] If $\mathbf{X}$ is a $2$-Hermitian matrix,  then  its eigenvalues taking the form $\lambda=\frac{\lambda_s+\lambda_d}{2}+\frac{\lambda_s-\lambda_d}{2\sqrt{\beta}}\mathrm{j}$,  where $\lambda_\nu\in \mathbb{R}$, for $\nu=s,d$,  are the eigenvalues of $\mathbf{X}_s$ and $\mathbf{X}_d$, respectively. Furthermore,  
equation \eqref{eigdes} takes the form
    \begin{equation*}
\mathbf{X}=\mathbf{U}\mathbf{\Lambda}\mathbf{U}^{\texttt{H}_{-1}^{1}}
\end{equation*}
\item[ \normalfont{Case }$\alpha>0$:]
If $\mathbf{X}$ is a $1$-Hermitian matrix,  then  its eigenvalues and eigenvectors are $(\alpha\beta)$-tessarines.
    Furthermore,   equation \eqref{eigdes} takes the form
    \begin{equation}\label{posher}
\mathbf{X}=\mathbf{U}\mathbf{\Lambda}\mathbf{U}^{\texttt{H}_{1}^{1}}
\end{equation}
\end{description}
\end{corollary}

\begin{proof}
 Next,  case $\alpha>0$ is demonstrated. The proof for case $\alpha<0$ follows analogously and ii has been omitted for brevity.

From Property 
\ref{proherm1}, matrices $\mathbf{A}_n$,  for $n=1,\ldots,4$,  given in \eqref{matricessAA},  are symmetric, and hence,   their eigenvalues and eigenvectors are real-valued.  Thus, from Remark \ref{foreig},  it follows that the eigenvalues and eigenvectors of $\mathbf{X}$ are $(\alpha\beta)$-tessarines.   Additionally,  matrices $\mathbf{V}_n$,  for $n=1,\ldots,4$,  given in \eqref{matV},  satisfy relation $\mathbf{V}_n^{-1}=\mathbf{V}_n^{\texttt{T}}$, for all $n=1,\ldots,4$. Finally, the result \eqref{posher} is  obtained from \eqref{matV}, Corollary \ref{injnj},  and Property    
\ref{proherm1}.
\end{proof}

\begin{theorem}
 Suppose that $\mathbf{X}\in \mathbb{T}_{\alpha\beta}^{p\times p}$ is $n$-Hermitian for $n=1,2$.  Then,  the following three conditions are equivalent:
\begin{enumerate}
\item $\mathbf{X}$ is positive  definite.
\item The eigenvalues of $\mathbf{X}$ are positive.
\item The real matrices $\mathbf{A}_n$, $n=1,\ldots,4$ (for $\alpha>0$), given in \eqref{matricessAA},  and the $(\alpha)$-complex matrices $\mathbf{X}_\nu$, $\nu=s,d$  (for  $\alpha<0$),  given in Proposition \ref{prpp2},   are positive definite.
\end{enumerate}
\end{theorem}

\begin{proof}
Next, case $\alpha>0$ is demonstrated. The proof for case $\alpha<0$ follows in an analogous manner and it has been omitted for brevity.

First,  it is demonstrated that point 1 implies point 2.  Consider an eigenvalue $\lambda$ of $\mathbf{X}$ and its corresponding eigenvector $\mathbf{u}$.
From the expressions \eqref{eigenvnv} and \eqref{matV},  and  Corollary \ref{Cor3},   it follows that  the associated $(\alpha\beta)$-tessarine of $\lambda$  
   is  $ \theta_1+\theta_2\mathrm{i}+ \theta_3\mathrm{j}+\theta_4 \mathrm{k} $,   $ \theta_n\in \mathbb{R}$, $n=1,\dots,4$,  
   and the  associated $(\alpha\beta)$-tessarine vector of $\mathbf{u}$  
   is $\mathbf{v}_1+\mathbf{v}_2\mathrm{i}+\mathbf{v}_3\mathrm{j}+\mathbf{v}_4\mathrm{k}  $.

 Now,  consider that $\mathbf{u}^{\texttt{H}_{1}^{1}}\mathbf{X}\mathbf{u}=\lambda\mathbf{u}^{\texttt{H}_{1}^{1}}\mathbf{u}$.  Applying Property \ref{proherm1},   $\lambda\mathbf{u}^{\texttt{H}_{1}^{1}}\mathbf{u}$ has associated $(\alpha\beta)$-tessarine:
 $\theta_1\mathbf{v}_1^\texttt{T}\mathbf{v}_1+ \theta_2\mathbf{v}_2^\texttt{T}\mathbf{v}_2\mathrm{i}+ \theta_3\mathbf{v}_3^\texttt{T}\mathbf{v}_3\mathrm{j}+\theta_4\mathbf{v}_4^\texttt{T}\mathbf{v}_4 \mathrm{k} $, 
and  taking into account  that $\mathbf{X}$ is positive definite,  it can be concluded that $\lambda$ is positive.

Next, it is demonstrated that point 2 implies point 1. From \eqref{posher},  the following expression holds:  $$\mathbf{x}^{\texttt{H}_{1}^{1}}\mathbf{X}\mathbf{x}=
\mathbf{x}^{\texttt{H}_{1}^{1}}\mathbf{U}\mathbf{\Lambda}\mathbf{U}^{\texttt{H}_{1}^{1}}\mathbf{x}=\mathbf{y}^{\texttt{H}_{1}^{1}}
\mathbf{\Lambda}\mathbf{y}=\sum_{i=1}^p\lambda_iy_i^2$$ 
By applying similar reasoning as before,  it can be concluded that  $\mathbf{x}^{\texttt{H}_{1}^{1}}\mathbf{X}\mathbf{x}$ is positive.

The proof of the equivalence of point 3 is straightforward. 

\end{proof}

\begin{corollary}
 Suppose that $\mathbf{X}\in \mathbb{T}_{\alpha\beta}^{p\times p}$ is positive  definite. Then, a unique positive definite square root matrix exists whose eigenvalues are the square roots of the eigenvalues of $\mathbf{X}$.
\end{corollary}

\subsubsection{Power Method}

The power method is a well-established iterative algorithm for eigenvalue computation,  widely used in real and complex domains. In this section,  the application of this method is extended to the  $(\alpha\beta)$-tessarine domain. 
The  behavior and  convergence properties of this method are analyzed  within this algebraic framework,  demonstrating that the algorithm converges to one of the eigenvalues of a given $(\alpha\beta)$-tessarine matrix $\mathbf{X}\in \mathbb{T}_{\alpha\beta}^{p\times p}$,  referred to as the dominant eigenvalue.  This extension provides a valuable tool for the optimal spectral representation of $\mathbf{X}$.

Consider a matrix $\mathbf{X}\in \mathbb{T}_{\alpha\beta}^{p\times p}$ that  is  $n$-Hermitian,  with $n=1,2$.  The following assumptions are also made:
\begin{itemize}
\item For $\alpha<0$,  the eigenvalues of  $\mathbf{X}_\nu$,  with $\nu=s,d$,  are ordered as follows:
\begin{equation*}\label{orhne}
|\lambda_{\nu1}|> |\lambda_{\nu2}|>|\lambda_{\nu 3}|>\ldots>|\lambda_{\nu p}|, \qquad \nu=s,d
\end{equation*}
 \item For $\alpha>0$,  the eigenvalues of the matrices  $\mathbf{A}_n$,  with $n=1,\ldots,4$,  defined  in \eqref{matricessAA},  are ordered as follows:
\begin{equation}\label{hnhnbvm}
|\theta_{n1}|> |\theta_{n2}|>|\theta_{n 3}|>\ldots>|\theta_{n p}|, \qquad n=1,\ldots,4
\end{equation}
\end{itemize}

Given $\mathbf{x} \in \mathbb{T}_{\alpha\beta}^p$,   define $|\mathbf{x}|_n=\left(\mathbf{x}^{\texttt{H}_{3-2n}^{1}}\mathbf{x}\right)^{1/2}$,  with $n=1$ for $\alpha>0$ and $n=2$ for $\alpha<0$.

\begin{property}\label{propmodu}
 Let $\mathbf{x}, \mathbf{y} \in \mathbb{T}_{\alpha\beta}^p$.  For $n=1,2$,  the following properties hold: 
\begin{enumerate}
\item $|\mathbf{x}|_n$ is an $(\alpha\beta)$-tessarine.
\item If $x \in \mathbb{T}_{\alpha\beta}$,   then   $|x\mathbf{y}|_n=|x|_n|\mathbf{y}|_n$.
\item $\displaystyle\left|\frac{\mathbf{x}}{|\mathbf{x}|_n}\right|_n=1.$
\item $|\det(\mathbf{P})|_n=1$,   with  $\det(\mathbf{P})$ given in Table \ref{tabdt}.
\end{enumerate}
\end{property}

Finally, consider the eigenvalue $\lambda_1$ of $\mathbf{X}\in \mathbb{T}_{\alpha\beta}^{p\times p}$,  obtained from eigenvalues $\lambda_{\nu1}$, $\nu=s,d$, for $\alpha<0$  (see \eqref{ealphane}),   and eigenvalues $\theta_{n1}$, $n=1,\ldots,4$, for case $\alpha>0$ (see \eqref{eigenvnv}). This eigenvalue is referred to as  the dominant eigenvalue,  and its corresponding eigenvector is denoted by $\mathbf{u}_1$.

The following result analyzes the convergence of the iterative sequence
$\mathbf{x}_m=\displaystyle\frac{\mathbf{X}\mathbf{x}_{m-1}}{|\mathbf{X}\mathbf{x}_{m-1}|_n}$, $m>0$,   with  $\mathbf{x}_0$ such that $|\mathbf{x}_0|_n=1$,  and its relation to $\lambda_1$.

\begin{theorem}\label{powerla}
Under the previous conditions,  the following property holds:
\begin{equation}\label{dfgrhgo}
\mathbf{x}_m^{\texttt{H}_{3-2n}^{1}}\mathbf{X}\mathbf{x}_m\rightarrow \lambda_1
\end{equation}
as $m$ tends to $\infty$,  for   $n=1,2$.
\end{theorem}

\begin{proof}
A proof is provided for case  $\alpha>0$.  The proof for case $\alpha<0$ follows analogously and it has been  omitted for brevity.

Let $\{\lambda_i\}_{i=2}^p$  be any set of $p-1$ eigenvalues of $\mathbf{X}$ such that,  together with $\lambda_1$,  they satisfy the conditions of Theorem \ref{prooo3}.  The corresponding eigenvectors are denoted by $\{\mathbf{u}_i\}_{i=1}^p$.   Consider any $\mathbf{x}_0\in \mathbb{T}_{\alpha\beta}^p $ such that $|\mathbf{x}_0|_1=1$. Then,
\begin{equation*}
\mathbf{x}_0=\sum_{i=1}^pa_i\mathbf{u}_i, \quad a_i \in \mathbb{T}_{\alpha\beta}, \quad i=1,\ldots,p
\end{equation*}

Since  $\mathbf{X}\mathbf{u}_i=\lambda_i\mathbf{u}_i$, it follows that
\begin{equation*}
\mathbf{X}^m\mathbf{x}_0=a_1\lambda_1^m\left(\mathbf{u}_1+\sum_{i=2}^p\frac{a_i}{a_1}
\left(\frac{\lambda_i}{\lambda_1}\right)^m \mathbf{u}_i\right)
\end{equation*}
Thus, $\mathbf{x}_m$ can be expressed as
$\mathbf{x}_m=C_{m_1}C_{m_2}$,  where $C_{m_1}=\frac{a_1\lambda_1^m}{|a_1\lambda_1^m|_1}$ and
\begin{equation*}
C_{m_2}=\frac{\mathbf{u}_1+\displaystyle\sum_{i=2}^p\frac{a_i}{a_1}
\left(\frac{\lambda_i}{\lambda_1}\right)^m \mathbf{u}_i}{\left| \mathbf{u}_1+\displaystyle\sum_{i=2}^p\frac{a_i}{a_1}
\left(\frac{\lambda_i}{\lambda_1}\right)^m \mathbf{u}_i\right|_1}
\end{equation*}

In accordance with Property \ref{plpprop1}, term  $\left(\frac{\lambda_i}{\lambda_1}\right)^m $
has the associated  $(\alpha\beta)$-tessarine
\begin{equation}\label{RR1}
\left(\frac{\tau_1(\theta_{1i})}{\theta_{11}}\right)^m+\left(\frac{\tau_2(\theta_{2i})}{\theta_{21}}\right)^m\mathrm{i}
\left(\frac{\tau_3(\theta_{3i})}{\theta_{31}}\right)^m\mathrm{j}+\left(\frac{\tau_4(\theta_{4i})}{\theta_{41}}\right)^m\mathrm{k}, \quad i=2,\ldots,p
\end{equation}

Since  $\theta_{n1}$,  for $n=1,\ldots,4$,  are the largest eigenvalues,  the expression \eqref{RR1} converges to zero as $m$ increases.  Consequently,  $C_{m_2}$ converges to  $\mathbf{u}_1$. Furthermore,  from Property \ref{propmodu}, sequence $C_{m_1}$ is bounded by $1$,   thereby achieving the result in \eqref{dfgrhgo}.
\end{proof}

\begin{remark}
For a sufficiently large $m$, it follows that  $\mathbf{x}_m\cong \pi\mathbf{u}_1$, where $\pi$  satisfies  $|\pi|_n=1$ and hence, $\mathbf{x}_m$ is an element of  the eigenspace associated to $\mathbf{u}_1$.
\end{remark}

\begin{remark}
An important consequence of  Theorem \ref{powerla} is that it provides a method for obtaining a rank-$k$ approximation for the eigenrepresentation of a positive definite matrix (see Section \ref{subnxm}).


\end{remark}

\subsection{SVD of an $(\alpha\beta)$-Tessarine Matrix}
This section addresses the SVD   of an $(\alpha\beta)$-tessarine matrix.
\begin{proposition}\label{ssvvdd2}
Consider $\mathbf{X}\in \mathbb{T}_{\alpha\beta}^{p\times q}$ is expressed in form  \eqref{repma} as $\mathbf{X}=\mathbf{X}_{s}\mathrm{w}_1+\mathbf{X}_{d}\mathrm{w}_2$.  Let $\mathbf{X}_s=\mathbf{U}_s\mathbf{\Sigma}_s\mathbf{V}_s^{\texttt{H}_{3-2n}^1}$ and $\mathbf{X}_d=\mathbf{U}_d\mathbf{\Sigma}_d\mathbf{V}_d^{\texttt{H}_{3-2n}^1}$ be  the SVDs  of $\mathbf{X}_s$ and $\mathbf{X}_d$ (with $n=1$ for $\alpha>0$ and $n=2$ for $\alpha<0$), respectively, computed from  Lemma \ref{SVDlema}. Then, $\mathbf{X}=\mathbf{U}\mathbf{\Sigma}\mathbf{V}^{\texttt{H}_{3-2n}^1}$  with the matrices
        \begin{equation}\label{mnhkyuierfg}
        \begin{split}
        \mathbf{U}&=\frac{\mathbf{U}_{s}+\mathbf{U}_{d}}{2}
        +\frac{\mathbf{U}_{s}-\mathbf{U}_{d}}{2\sqrt{\beta}}\mathrm{j}\\
          \mathbf{\Sigma}&=\frac{\mathbf{\Sigma}_{s}+\mathbf{\Sigma}_{d}}{2}
        +\frac{\mathbf{\Sigma}_{s}-\mathbf{\Sigma}_{d}}{2\sqrt{\beta}}\mathrm{j}\\
          \mathbf{V}&=\frac{\mathbf{V}_{s}+\mathbf{V}_{d}}{2}
        +\frac{\mathbf{V}_{s}-\mathbf{V}_{d}}{2\sqrt{\beta}}\mathrm{j}
        \end{split}
        \end{equation}
\end{proposition}

\begin{remark}

If $\alpha<0$,  the singular values, i.e, the non-zero elements of $\mathbf{\Sigma}$ are $(\alpha\beta)$-tessarines numbers in which components $\mathrm{i}$ and $\mathrm{k}$ vanish.   Moreover, if $r_s$ and $r_d$ denote the numbers of singular values associated with the SVDs of $\mathbf{X}_s$ and $\mathbf{X}_d$, respectively,  then  there are $r^2$ possible SVDs of $\mathbf{X}$, with $r=\max(r_s,r_d)$. 

However, if $\alpha>0$,  the singular values of $\mathbf{\Sigma}$  are $(\alpha\beta)$-tessarines.   Moreover,  if $r_n$,  for $n=1,\ldots,4$,  denote the numbers of singular values associated with the SVDs of the matrices $\mathbf{A}_n$,   $n=1,\ldots,4$, respectively,   defined in \eqref{matricessAA},   then  there are $r^4$ possible SVDs  of $\mathbf{X}$,   with $r=\max(r_1,r_2,r_3,r_4)$.
\end{remark}

\subsubsection{Rank-$k$ Approximation}\label{subnxm}
The goal of this section is to determine the best rank-$k$ approximation of an $(\alpha\beta)$-tessarine matrix $\mathbf{X}$.  Initially,  the best approximation among the different existing SVD representations of  $\mathbf{X}$  is analyzed.  Then,  the general case of any rank-$k$ matrix is addressed.

For simplicity and clarity of exposition,  the case in which $\alpha>0$ is analyzed first.  Subsequently,  a brief analysis is provided on how the results derived for this case can be extrapolated to the case in which $\alpha<0$.

Let 
\begin{equation}\label{skingvalk}
\sigma_{n1}\geq \sigma_{n2}\geq \ldots \geq \sigma_{nr_n}, \qquad n=1,\ldots,4
\end{equation}
denote  the set of singular values of matrices $\mathbf{A}_n$ defined in   \eqref{matricessAA},  and consider their  SVD representations:
 \begin{equation}\label{fgsvd}
 \mathbf{A}_n=\mathbf{U}_n\mathbf{\Sigma}_n\mathbf{V}_n^{\texttt{T}}
   \end{equation}
where    $\sigma_{ni}$ is  the $(i,i)$th  element of $\mathbf{\Sigma}_n$.  Consider the SVD representation of $\mathbf{X}$ obtained from the decomposition  $\mathbf{X}=\mathbf{U}\mathbf{\Sigma}\mathbf{V}^{\texttt{H}_{1}^1}$.    It can be clearly demonstrated that  $\mathbf{X}=\displaystyle\sum_{i=1}^r\sigma_i\mathbf{u}(i)\mathbf{v}^{\texttt{H}_{1}^{1}}(i)$, 
where  $r=\max(r_1,r_2,r_3,r_4)$,  $\sigma_i$ is the $(i,i)$th element of $\boldsymbol{\Sigma}$,  and $\mathbf{u}(i)$ and $\mathbf{v}(i)$ are the $i$th  columns of $\mathbf{U}$ and $\mathbf{V}$, respectively.  
 
The following definition is introduced:
 \begin{equation}\label{xkxk}
 \mathbf{X}_k=\displaystyle\sum_{i=1}^k\sigma_i\mathbf{u}(i)\mathbf{v}^{\texttt{H}_{1}^{1}}(i), \qquad k\leq r
\end{equation}
Moreover,  consider four permutations $\{\tau_n(\sigma_{ni})\}_{i=1}^{r_n}$,  for $n=1,\dots,4$,  which yield a  new SVD of  $\mathbf{X}$ given by  $\mathbf{X}=\check{\mathbf{U}}\check{\mathbf{\Sigma}}\check{\mathbf{V}}^{\texttt{H}_{1}^1}$.  Based on this decomposition,  the following expression is defined: 
  $\check{\mathbf{X}}_k=\displaystyle\sum_{i=1}^k\check{\sigma}_i\check{\mathbf{u}}(i)\check{\mathbf{v}}^{\texttt{H}_{1}^{1}}(i)$,   where $\check{\sigma}_i$ represents the $(i,i)$th element of $\check{\boldsymbol{\Sigma}}$,  and $\check{\mathbf{u}}(i)$ and $\check{\mathbf{v}}(i)$ correspond to the $i$th columns of $\check{\mathbf{U}}$ and $\check{\mathbf{V}}$, respectively.
\begin{theorem}\label{trrrt}
Let $\mathbf{X}\in \mathbb{T}_{\alpha\beta}^{p\times q} $.  The following properties hold:
\begin{enumerate}
\item $<\mathbf{X},\mathbf{X}>_1=\displaystyle\sum_{i=1}^r\sigma_i^2$.
\item  $<\mathbf{X}-\mathbf{X}_k,\mathbf{X}-\mathbf{X}_k>_1=\displaystyle\sum_{i=k+1}^r\sigma_i^2$.
\item  $ <\mathbf{X}-\mathbf{X}_k,\mathbf{X}-\mathbf{X}_k>_1\preceq <\mathbf{X}-\check{\mathbf{X}}_k,\mathbf{X}-\check{\mathbf{X}}_k>_1$.
    \item $\|\mathbf{X}-\mathbf{X}_k\|_1\leq\|\mathbf{X}-\check{\mathbf{X}}_k\|_1$.
\end{enumerate}
where  $<\cdot,\cdot>_1$ and $\|\cdot\|_1$ are defined  in \eqref{producttoma} and \eqref{normim}, respectively.
\end{theorem}

\begin{proof}
To demonstrate points 1 and 2, first, following equations need to be established:
 \begin{equation}\label{ortgij}
 \mathbf{U}^{\texttt{H}_{1}^{1}}\mathbf{U}=\mathbf{I}_p, \qquad \mathbf{V}^{\texttt{H}_{1}^{1}}\mathbf{V}=\mathbf{I}_q
 \end{equation}

 From Lemma \ref{SVDlema} and Proposition \ref{ssvvdd2}, the associated  $(\alpha\beta)$-tessarines matrices of  $\mathbf{U}$ and $\mathbf{V}$ are directly given by $\tilde{\mathbf{U}}=\mathbf{U}_1+\mathbf{U}_2\mathrm{i}+\mathbf{U}_3\mathrm{j}+\mathbf{U}_4\mathrm{k}$ and  $\tilde{\mathbf{V}}=\mathbf{V}_1+\mathbf{V}_2\mathrm{i}+\mathbf{V}_3\mathrm{j}+\mathbf{V}_4\mathrm{k}$, respectively,
where $\mathbf{U}_n$ and $\mathbf{V}_n$,  for $n=1,\ldots,4$,   are defined in \eqref{fgsvd}. Then, by applying Property \ref{proherm1},  the relations in  \eqref{ortgij} are obtained.  Consequently,  the properties stated in points 1 and 2 are devised.

Next, suppose without loss of generality that $r=r_1$. In this situation,  every singular value $\sigma_i$ has the associated $(\alpha\beta)$-tessarine $\sigma_{1i}+\sigma_{2i}\mathrm{i}+\sigma_{3i}\mathrm{j}+\sigma_{4i}\mathrm{k}$,  with $\sigma_{ni}=0$,  for $i>r_n$  and $n=2,3,4$. From Property \ref{plpprop1},  $\sigma_i^2$  has  the associated $(\alpha\beta)$-tessarine  $\sigma_{1i}^2+\sigma_{2i}^2\mathrm{i}+\sigma_{3i}^2\mathrm{j}+\sigma_{4i}^2\mathrm{k}$. Thus,  using equation \eqref{skingvalk},  point 3 follows.

Finally,  point 4 is a direct consequence of point 3 and Property \ref{realnrrr}.
\end{proof}


Hereafter,   let $\{\sigma_i(\mathbf{X})\}_{i=1}^r$ denote the nonincreasingly ordered singular values of a real matrix $\mathbf{X} \in \mathbb{R}^{p\times q}$,  and the set of singular values obtained via \eqref{skingvalk} for $(\alpha\beta)$-tessarine matrix $\mathbf{X}\in \mathbb{T}_{\alpha\beta}^{p\times q}$.

The general case of the rank-$k$ approximation of a matrix is addressed below.
\begin{property}\label{rangoji}
Let $\mathbf{X} \in \mathbb{T}_{\alpha\beta}^{p\times q}$  denote its  associated $(\alpha\beta)$-tessarine
matrix by $\tilde{\mathbf{X}}=\tilde{\mathbf{X}}_1+\tilde{\mathbf{X}}_2\mathrm{i}+\tilde{\mathbf{X}}_3\mathrm{j}+\tilde{\mathbf{X}}_4\mathrm{k}$,  and let $\vartheta_i$ represent  the associated $(\alpha\beta)$-tessarine of the singular value $\sigma_i(\mathbf{X})$.  Then, the following properties hold:
\begin{enumerate}
\item \begin{equation}\label{derfdg}
\vartheta_i=\sigma_i(\tilde{\mathbf{X}}_1)+\sigma_i(\tilde{\mathbf{X}}_2)\mathrm{i}+\sigma_i(\tilde{\mathbf{X}}_3)\mathrm{j}+\sigma_i(\tilde{\mathbf{X}}_4)\mathrm{k}
\end{equation}
\item $\rank(\mathbf{X})$  coincides with the maximum number of positive singular values of $\tilde{\mathbf{X}}_n$, for $n=1,\ldots,4$.
\item $\mathbf{X}$ and $\mathbf{X}^{\texttt{H}_{1}^{1}}$ share the same set of singular values.
\end{enumerate}
\end{property}

\begin{theorem}\label{zxcvbn}
Let $\mathbf{X}, \mathbf{Y} \in \mathbb{T}_{\alpha\beta}^{p\times q}$,  and $r=\min(p,q)$.  The following results are verified:
\begin{enumerate}
\item
\begin{equation}\label{klpoiop}
 <\mathbf{X}^{\texttt{H}_1^1},\mathbf{Y}^{\texttt{H}_1^1}>_1\preceq \sum_{i=1}^r \sigma_i(\mathbf{X})\sigma_i(\mathbf{Y})
 \end{equation}
    \item
    \begin{equation}\label{xxxxkl}
    <\mathbf{X}-\mathbf{Y},\mathbf{X}-\mathbf{Y}>_1
    \succeq \displaystyle \sum_{i=1}^r (\sigma_i(\mathbf{X})-\sigma_i(\mathbf{Y}))^2
    \end{equation}
\end{enumerate}
\end{theorem}
\begin{proof}
Let the associated $(\alpha\beta)$-tessarine of $ <\mathbf{X}^{\texttt{H}_1^1},\mathbf{Y}^{\texttt{H}_1^1}>_1=\trace(\mathbf{X}\mathbf{Y}^{\texttt{H}_{1}^{1}})$ be denoted by $y$,  and let the associated $(\alpha\beta)$-tessarine matrices of $\mathbf{X}\mathbf{Y}^{\texttt{H}_{1}^{1}}$, $\mathbf{X}$,  and $\mathbf{Y}^{\texttt{H}_{1}^{1}}$ be denoted  by $\mathbf{E}$, $\tilde{\mathbf{X}}=\tilde{\mathbf{X}}_1+\tilde{\mathbf{X}}_2\mathrm{i}+\tilde{\mathbf{X}}_3\mathrm{j}+\tilde{\mathbf{X}}_4\mathrm{k}$ and $\mathbf{T}={\mathbf{T}}_1+{\mathbf{T}}_2\mathrm{i}+{\mathbf{T}}_3\mathrm{j}+{\mathbf{T}}_4\mathrm{k}$, respectively.

From Property \ref{proherm1}, the following equation holds:
\begin{equation*}
y=\trace(\mathbf{E})=\trace(\tilde{\mathbf{X}}_1{\mathbf{T}}_1)+\trace(\tilde{\mathbf{X}}_2{\mathbf{T}}_2)\mathrm{i}+\trace(\tilde{\mathbf{X}}_3{\mathbf{T}}_3)\mathrm{j}
+\trace(\tilde{\mathbf{X}}_4{\mathbf{T}}_4)\mathrm{k}
\end{equation*}
Then,  by  applying the von Neumann inequality \cite{Horn},   expression \eqref{derfdg},  and Properties \ref{rangoji} and  \ref{proherm1},  the relation \eqref{klpoiop} is devised.

Next, from Property \ref{rodproxpro}, the following expression holds:
\begin{equation*}
<\mathbf{X}-\mathbf{Y},\mathbf{X}-\mathbf{Y}>_1=<\mathbf{X},\mathbf{X}>_1+<\mathbf{Y},\mathbf{Y}>_1-
2<\mathbf{X}^{\texttt{H}_1^1},\mathbf{Y}^{\texttt{H}_1^1}>_1
\end{equation*}

 Finally,  by considering Theorem \ref{trrrt} and \eqref{klpoiop}, relation  \eqref{xxxxkl} is proven.

\end{proof}

\begin{theorem} \label{ssaazzxxcv}
Let $\mathbf{X}, \mathbf{Y} \in \mathbb{T}_{\alpha\beta}^{p\times q}$,   and  $r=\min(p,q)$.  Suppose that $\rank(\mathbf{Y})=k$, $1\leq k<r$.  The following results are verified:
\begin{enumerate}
\item  $ <\mathbf{X}-\mathbf{X}_k,\mathbf{X}-\mathbf{X}_k>_1\preceq <\mathbf{X}-\mathbf{Y},\mathbf{X}-\mathbf{Y}>_1$.
    \item $\|\mathbf{X}-\mathbf{X}_k\|_1\leq\|\mathbf{X}-\mathbf{Y}\|_1$.
\end{enumerate}
with $\mathbf{X}_k$ given in \eqref{xkxk}.
\end{theorem}
\begin{proof}
From Theorem \ref{zxcvbn} and Property \ref{rangoji}, the following relations hold:
\begin{equation*}
\begin{split}
<\mathbf{X}-\mathbf{Y},\mathbf{X}-\mathbf{Y}>_1
    \succeq  &\sum_{i=1}^r (\sigma_i(\mathbf{X})-\sigma_i(\mathbf{Y}))^2\\
    =&\sum_{i=1}^k(\sigma_i(\mathbf{X})-\sigma_i(\mathbf{Y}))^2+\sum_{i=k+1}^r \sigma_i^2(\mathbf{X})\succeq
    \sum_{i=k+1}^r \sigma_i^2(\mathbf{X})
 \end{split}
    \end{equation*}
By applying Theorem \ref {trrrt},  point 1  is demonstrated.   Finally,   point 2 follows directly from point 1 and Property \ref{realnrrr}. 
\end{proof}

\begin{corollary}\label{corolid}
Let  $\mathbf{X}, \mathbf{Y} \in \mathbb{T}_{\alpha\beta}^{p\times p}$. Suppose that $\mathbf{X}$ is positive definite  and $\rank(\mathbf{Y})=k$, $1\leq k<p$.   Let $\{\lambda_i\}_{i=1}^p$ denote the set of eigenvalues of $\mathbf{X}$,  ordered in accordance with  \eqref{hnhnbvm},  and define $\mathbf{X}_k=\sum_{i=1}^k\lambda_i\mathbf{u}(i)\mathbf{u}^{\texttt{H}_{1}^{1}}(i)$, where $\mathbf{u}(i)$ is the corresponding eigenvector of $\lambda_i$.  The following results are verified:
\begin{enumerate}
\item  $ <\mathbf{X}-\mathbf{X}_k,\mathbf{X}-\mathbf{X}_k>_1=\displaystyle\sum_{i=k+1}^p\lambda_i\preceq <\mathbf{X}-\mathbf{Y},\mathbf{X}-\mathbf{Y}>_1$.
    \item $\|\mathbf{X}-\mathbf{X}_k\|_1\leq\|\mathbf{X}-\mathbf{Y}\|_1$.
\end{enumerate}
\end{corollary}

Next,   the previously obtained results are adapted to case $\alpha<0$.  For this purpose, consider
\begin{equation}\label{fgfgd}
\sigma_{\nu 1}\geq \sigma_{\nu 2}\geq \ldots \geq \sigma_{\nu r_{\nu}}, \qquad \nu=s,d
\end{equation}
 as the set of singular values of matrices $\mathbf{X}_s$ and $\mathbf{X}_d$ given in \eqref{repma}.
Define  the SVD of $\mathbf{X}$ associated to \eqref{fgfgd} as $$\mathbf{X}=\displaystyle\sum_{i=1}^r\sigma_i\mathbf{u}(i)\mathbf{v}^{\texttt{H}_{-1}^{1}}(i)$$  with $r=\max(r_s,r_d)$.

Theorem \ref{trrrt} and points 2 and 3 of Property \ref{rangoji} can be clearly proven in the case where $\alpha<0$.  The corresponding version of Theorem \ref{zxcvbn} is presented below.

\begin{theorem}\label{zxcxcvvbn}
Let  $\mathbf{X}, \mathbf{Y} \in \mathbb{T}_{\alpha\beta}^{p\times q}$,  and $r=\min(p,q)$.   Then,   the following inequalities hold:
\begin{enumerate}
\item
\begin{equation}\label{klggpoiop}
\Re\{ <\mathbf{X}^{\texttt{H}_{-1}^1},\mathbf{Y}^{\texttt{H}_{-1}^1}>_2\}\leq\Re\left\{\sum_{i=1}^r \sigma_i(\mathbf{X})\sigma_i(\mathbf{Y})\right\}
 \end{equation}
    \item
    \begin{equation}\label{gggghhhhj}
    \|\mathbf{X}-\mathbf{Y}\|_2^2
    \geq \Re\left\{\displaystyle \sum_{i=1}^r (\sigma_i(\mathbf{X})-\sigma_i(\mathbf{Y}))^2\right\}
    \end{equation}
\end{enumerate}
\end{theorem}
\begin{proof}
Inequality \eqref{klggpoiop} can be verified for matrices taking the form  $\mathbf{X}=\mathbf{A}_1+\mathbf{B}_1\mathrm{i}$  and $\mathbf{Y}=\mathbf{A}_2+\mathbf{B}_2\mathrm{i}$ by applying  the von Neumann inequality for complex matrices.   

For the general case,  using \eqref{repma},  it follows that:
\begin{multline*}
\Re\{ <\mathbf{X}^{\texttt{H}_{-1}^1},\mathbf{Y}^{\texttt{H}_{-1}^1}>_2\}=\Re\{\trace(\mathbf{X}_s\mathbf{Y}_s^{\texttt{H}_{-1}^1}\mathrm{w}_1+
\mathbf{X}_d\mathbf{Y}_d^{\texttt{H}_{-1}^1}\mathrm{w}_2)\}\\ =\frac{1}{2}\Re\{\trace(\mathbf{X}_s\mathbf{Y}_s^{\texttt{H}_{-1}^1}+
\mathbf{X}_d\mathbf{Y}_d^{\texttt{H}_{-1}^1})\}
\end{multline*}
Therefore,  by considering the above expression and applying \eqref{mnhkyuierfg},  inequality \eqref{klggpoiop} is established. 

The proof of \eqref{gggghhhhj} follows similarly from \eqref{xxxxkl}.
\end{proof}

A version of point 2 of Theorem \ref{ssaazzxxcv} can be derived by applying Property \ref{rodproxpro} and Theorem \ref{zxcxcvvbn}.  Consequently,  a result similar to that of point 2 of Corollary \ref{corolid} is obtained.

\subsection{Pseudoinverse  of an $(\alpha\beta)$-Tessarine Matrix}
To conclude this section, the concept of the $(\alpha\beta)$-tessarine pseudoinverse matrix is introduced.  This concept will be useful in the next section.
\begin{defi}\label{propseudoiy}
Let  $\mathbf{X}\in \mathbb{T}_{\alpha\beta}^{p\times q}$.
The $n$-pseudoinverse of  $\mathbf{X}$, where $n=1$ for $\alpha>0$ and $n=2$ for $\alpha<0$,  denoted as   $\mathbf{P}$,   is an $(\alpha\beta)$-tessarine matrix  satisfying the following four conditions:
\begin{enumerate}
\item $\mathbf{X}\mathbf{P}\mathbf{X}=\mathbf{X}$.
\item $\mathbf{P}\mathbf{X}\mathbf{P}=\mathbf{P}$.
\item $(\mathbf{X}\mathbf{P})^{\texttt{H}_{3-2n}^{1}}=\mathbf{X}\mathbf{P}$.
\item $(\mathbf{P}\mathbf{X})^{\texttt{H}_{3-2n}^{1}}=\mathbf{P}\mathbf{X}$.
\end{enumerate}
\end{defi}

The following result establishes the uniqueness of P and provides an explicit formula based on the  SVD of  $\mathbf{X}$.
\begin{theorem}
 Let $\mathbf{P}$ be the $n$-pseudoinverse of    $\mathbf{X}\in \mathbb{T}_{\alpha\beta}^{p\times q}$. Then,  the following properties hold: 
\begin{enumerate}
\item $\mathbf{P}$ is unique.
\item Consider matrices $\mathbf{U}$, $\mathbf{V}$, $\mathbf{\Sigma}_{s}$, and $\mathbf{\Sigma}_{d}$ given in Proposition \ref{ssvvdd2}.  Suppose that for $\alpha>0$,  it follows that
    \begin{equation*}
    \begin{split}
     \mathbf{\Sigma}_s&=\frac{\mathbf{\Sigma}_{s1}+\mathbf{\Sigma}_{d1}}{2}
        +\frac{\mathbf{\Sigma}_{s1}-\mathbf{\Sigma}_{d1}}{2\sqrt{\alpha}}\mathrm{i}\\
      \mathbf{\Sigma}_d&=\frac{\mathbf{\Sigma}_{s2}+\mathbf{\Sigma}_{d2}}{2}
        +\frac{\mathbf{\Sigma}_{s2}-\mathbf{\Sigma}_{d2}}{2\sqrt{\alpha}}\mathrm{i}
           \end{split}
           \end{equation*}
as obtained via Lemma \ref{SVDlema}. Then, the $n$-pseudoinverse of $\mathbf{X}$ is given by
\begin{equation}\label{pseudhiu}
\mathbf{P}=\mathbf{V}\mathbf{Q}\mathbf{U}^{\texttt{H}_{3-2n}^{1}}
\end{equation}
with  $\mathbf{Q}=\mathbf{P}_{s}\mathrm{w}_1+\mathbf{P}_{d}\mathrm{w}_2$,   where   
\begin{itemize}
\item For $\alpha<0$,  $\mathbf{P}_{s}$ and $\mathbf{P}_{d}$ are the pseudoinverses of the real matrices $\mathbf{\Sigma}_s$ and $\mathbf{\Sigma}_d$,  respectively.  

\item For $\alpha>0$,  
    \begin{equation*}
    \begin{split}
     \mathbf{P}_s&=\frac{\mathbf{P}_{s1}+\mathbf{P}_{d1}}{2}
        +\frac{\mathbf{P}_{s1}-\mathbf{P}_{d1}}{2\sqrt{\alpha}}\mathrm{i}\\
          \mathbf{P}_d&=\frac{\mathbf{P}_{s2}+\mathbf{P}_{d2}}{2}
        +\frac{\mathbf{P}_{s2}-\mathbf{P}_{d2}}{2\sqrt{\alpha}}\mathrm{i}
           \end{split}
           \end{equation*}
where $\mathbf{P}_{\nu n}$ denotes the pseudoinverses of the real matrices $\mathbf{\Sigma}_{\nu n}$,  for $\nu=s,d$  and  $n=1,2$.
\end{itemize}
\end{enumerate}
\end{theorem}

\begin{proof}
Suppose that two $n$-pseudoinverses,  $\mathbf{P}_1$ and $\mathbf{P}_2$, exist. By applying the conditions given in the definition, the following identity is obtained:
\begin{equation*}
\begin{split}
\mathbf{X}\mathbf{P}_1 & =\mathbf{X}\mathbf{P}_2\mathbf{X}\mathbf{P}_1=
(\mathbf{X}\mathbf{P}_2)^{\texttt{H}_{3-2n}^{1}}(\mathbf{X}\mathbf{P}_1)^{\texttt{H}_{3-2n}^{1}}=
\mathbf{P}_2^{\texttt{H}_{3-2n}^{1}} (\mathbf{X}\mathbf{P}_1\mathbf{X})^{\texttt{H}_{3-2n}^{1}}\\ & =
(\mathbf{X}\mathbf{P})_2^{\texttt{H}_{3-2n}^{1}}=\mathbf{X}\mathbf{P}_2
\end{split}
\end{equation*}

Similarly,  it follows that $\mathbf{P}_1\mathbf{X}=\mathbf{P}_2\mathbf{X}$.  Consequently,
\begin{equation*}
\mathbf{P}_1=\mathbf{P}_1\mathbf{X}\mathbf{P}_1=\mathbf{P}_1\mathbf{X}\mathbf{P}_2=\mathbf{P}_2\mathbf{X}\mathbf{P}_2=\mathbf{P}_2
\end{equation*} 

Finally, the expression in \eqref{pseudhiu} can be clearly verified as satisfying the four conditions stated in Definition \ref{propseudoiy}, thereby completing the proof.

\end{proof}

\section{The $(\alpha\beta)$-Tessarine Least Squares Problem}\label{sectiolineal}
In this section,   the least squares problem for $(\alpha\beta)$-tessarine matrices is addressed.
\begin{defi}\label{closedsett2}
Consider an $(\alpha\beta)$-tessarine matrix $\mathbf{X}\in\mathbb{T}_{\alpha\beta}^{p\times q}$. The closed linear subspace $\mathcal{D}$ corresponding to $\mathbf{X}$ is defined as the set of elements of the form $\sum_{i=1}^ q\phi_{i} \mathbf{x}(i)$, where $\phi_{i}\in \mathbb{T}_{\alpha\beta}$ are  deterministic numbers, and $ \mathbf{x}(i)$ denotes the $i$th column of $\mathbf{X}$.
\end{defi}

For any  $\mathbf{x},  \mathbf{y} \in \mathbb{T}_{\alpha\beta}^{p}$,  the following two distance measures can be derived from \eqref{normim}:
\begin{equation}\label{distdeter}
d_{n}(\mathbf{x},\mathbf{y})=\|\mathbf{x}-\mathbf{y}\|_{n},\qquad n=1,2
\end{equation}

Let  $\mathbf{y} \in \mathbb{T}_{\alpha\beta}^{p}$  be a deterministic vector.  The objective is
 to find a vector $\hat{\mathbf{y}}_{n}\in \mathcal{D}$ such that minimizes
 \begin{equation*}
 d_{n}(\mathbf{y},\mathbf{x}), \qquad \mathbf{x}\in\mathcal{D}
 \end{equation*}
where $n=1$ for  $\alpha>0$ and $n=2$ for  $\alpha<0$.

Vector $\hat{\mathbf{y}}_{n}$ is defined as the projection of $\mathbf{y} $ onto  the set $\mathcal{D}$ with respect to the distance measure given in \eqref{distdeter}.   In general,  in a metric space,   neither the existence nor the uniqueness of the projection is guaranteed.   However, following the methodology presented in \cite{Navarro_Segre}, the following result ensures both properties and provides a method for obtaining the projection.

\begin{theorem}\label{Teoremini}
Given $\mathbf{y} \in \mathbb{T}_{\alpha\beta}^{p}$,  and $\mathcal{D}$ the closed linear subspace associated to  $\mathbf{X}\in\mathbb{T}_{\alpha\beta}^{p\times q}$,   the following statements hold:
\begin{enumerate}
\item There exists a unique projection of $\mathbf{y} $, $\hat{\mathbf{y}}_{n}$, onto the set $\mathcal{D}$.
\item $\hat{\mathbf{y}}_{n}$ is the projection of $\mathbf{y} $ onto the set $\mathcal{D}$  if,  and only if,  $<\mathbf{x}, \mathbf{y} -\hat{\mathbf{y}}_{n}>_n=0$, $\forall \mathbf{x} \in \mathcal{D}$.
\end{enumerate}
\end{theorem}

\begin{corollary}
Given $\hat{\mathbf{y}}_{n}$ the projection of $\mathbf{y} \in \mathbb{T}_{\alpha\beta}^{p}$   onto the set  $\mathcal{D}$,  it follows that
\begin{equation*}
\hat{\mathbf{y}}_{n}=\mathbf{X}\mathbf{h}_n
\end{equation*}
with $\mathbf{X}\in\mathbb{T}_{\alpha\beta}^{p\times q}$ and 
\begin{equation}\label{parminv}
\begin{split}
\mathbf{h_n}&=(\mathbf{X}^{\texttt{H}_{3-2n}^{1}}\mathbf{X})^{-1}\mathbf{X}^{\texttt{H}_{3-2n}^{1}}\mathbf{y}\
\end{split}
\end{equation}
Furthermore,  the associated error $\varepsilon_{n}=\|\mathbf{y}-\hat{\mathbf{y}}_n\|_n^2$ is given by 
\begin{equation*}\label{parminv2}
\begin{split}
 \varepsilon_n &=\Re\left\{\mathbf{y}^{\texttt{H}_{3-2n}^{1}}\mathbf{y}
-\mathbf{y}^{\texttt{H}_{3-2n}^{1}}\mathbf{X}\mathbf{h}_n\right\}
\end{split}
\end{equation*}
\end{corollary}

In least squares problems,  ill-conditioned matrices lead to numerical instability,  making solutions highly sensitive to small data perturbations. To address this issue, the pseudoinverse provides a robust alternative to standard matrix inversion, ensuring  the existence of a least squares solution  even in  challenging scenarios.  
In this context,   a result for ill-conditioned matrices based on the $n$-pseudoinverse is presented in the $({\alpha\beta})$-tessarine domain. Specifically, by applying reasoning similar to that employed in the proof of Theorem 5 in \cite{Kosal1}, the following result is established.

\begin{corollary}  Vector $\mathbf{h}_n$ given in \eqref{parminv}, can be computed as follows:
\begin{equation}\label{leastpseuuu}
\mathbf{h_n}=\mathbf{P}\mathbf{y}
\end{equation}
with $\mathbf{y} \in \mathbb{T}_{\alpha\beta}^{p}$,  and where $\mathbf{P}$ is the $n$-pseudoinverse given in \eqref{pseudhiu}.
\end{corollary}

\subsection{Sequential Algorithm for Toeplitz $(\alpha\beta)$-Tessarine Matrices}

This section introduces a sequential algorithm for solving least squares problems involving $(\alpha\beta)$-tessarine matrices with an $n$-Hermitian Toeplitz structure (for further details, see \citep{Brockwell}). These systems are of interest because they are important in fields such as signal processing, image restoration, and neural network decoding \cite{Huang}.

Let  $\mathbf{X}\in\mathbb{T}_{\alpha\beta}^{p\times q}$ be a matrix such that
\begin{equation}\label{matrixT}
\mathbf{T}^{(n)}=\mathbf{X}^{\texttt{H}_{3-2n}^{1}}\mathbf{X}
\end{equation}
is an $n$-Hermitian Toeplitz matrix (see Definition \ref{toeplitz}),  where $n=1$ for $\alpha>0$ and $n=2$ for $\alpha<0$.  Define matrices $\mathbf{X}_1(i)=\left[\mathbf{x}(i),\ldots,\mathbf{x}(1)\right] $ and $\mathbf{X}_2(i)=\left[\mathbf{x}(2),\ldots,\mathbf{x}(i)\right]$,  with their  respective closed linear subspaces denoted by  $\mathcal{D}_1(i)$ and  $\mathcal{D}_2(i)$.

Two least squares problems within the $(\alpha\beta)$-tessarine domain are considered: the one-step updating problem and the fixed-point updating  problem.  In the first case,  the goal is to find  $\hat{\mathbf{x}}_{n}(i+1,i)\in \mathcal{D}_1(i)$ that minimizes
 \begin{equation*}
 d_{n}(\mathbf{x}(i+1),\mathbf{x}), \qquad \mathbf{x}\in\mathcal{D}_1(i)
 \end{equation*}
In the second case,  the objective is to find  $\hat{\mathbf{x}}_{n}(1,i)\in \mathcal{D}_2(i)$ that minimizes
 \begin{equation*}
 d_{n}(\mathbf{x}(1),\mathbf{x}), \qquad \mathbf{x}\in\mathcal{D}_2(i)
 \end{equation*}

Vectors $\hat{\mathbf{x}}_{n}(i+1,i)$ and $\hat{\mathbf{x}}_{n}(1,i)$ are referred to as  the one-stage and fixed-point projections,  respectively.   Assume that $\hat{\mathbf{x}}_{n}(i+1,i)= \mathbf{X}_1(i)\mathbf{f}_n(i)$ and $\hat{\mathbf{x}}_{n}(1,i)= \mathbf{X}_2(i)\mathbf{g}_n(i)$. Their associated errors  are  defined  by $ \varepsilon_{1n}(i)=\|\mathbf{x}(i+1)-\hat{\mathbf{x}}_{n}(i+1,i)\|_n^2$ and $  \varepsilon_{2n}(i)=\|\mathbf{x}(1)-\hat{\mathbf{x}}_{n}(1,i)\|_n^2$, respectively.  Furthermore,  define
\begin{equation*}
\begin{split}
\pi_n(i)&=\mathbf{T}_{1+i,1}^{(n)}, \qquad i\geq 0 \\
\pi_n(i)&= (\pi_n(-i))^{(3-2n,1)}, \qquad i<0
\end{split}
\end{equation*}

Due to the unique structure of $\mathbf{T}^{(n)}$,  the following proposition holds.
\begin{proposition}\label{proppest}
 For  $i\geq 2$,   the following relations hold:  
\begin{enumerate} 
 \item  $\mathbf{g}_n(i)=\mathbf{f}_n^{(3-2n,1)}(i-1)$.
 \item  $ \varepsilon_{2n}(i)= \varepsilon_{1n}(i-1)$.
 \end{enumerate}
\end{proposition}

To compute  $\mathbf{f}_n(i)$ and $\varepsilon_{1n}(i)$ efficiently as  $i$ increases,   a sequential algorithm is introduced.

\begin{theorem}\label{secuencialalgort}
Let  $ \varepsilon_{1n}(0)=\pi_n(0)$,   and define $\boldsymbol{\pi}_n(i)=[\pi_n(1-i),\ldots,\pi_n(-1)]$.  Suppose that $\mathbf{f}_n(i)=\left[
                                  \begin{array}{c}
                                  \boldsymbol{\lambda}_n(i)\\
                                    \lambda_n(i) \\
\end{array}
                                \right]$. Then,  the coefficients  $ \lambda_n(i)$  and $\boldsymbol{\lambda}_n(i)$  are given by                 
                               \begin{equation}\label{coeficn}
                                \begin{split}
                                \lambda_n(1)&=\pi_n(-1)\pi_n^{-1}(0)\\
                                 \lambda_n(i)&=\left[ \pi_n(-i)-\boldsymbol{\pi}_n(i)\mathbf{f}_n(i-1)\right]
                               \varepsilon_{1n}^{-1}(i-1), \qquad i\geq 2
                                \end{split}
                                \end{equation}
and 
    \begin{equation}\label{restodecoef}
  \begin{split}
  \boldsymbol{\lambda}_n(1)&=0\\
   \boldsymbol{\lambda}_n(i)&=  \mathbf{f}_n(i-1)-  \lambda_n(i) \boldsymbol{\rho}_n(i), \qquad i\geq 2
   \end{split}
   \end{equation}
where   $\boldsymbol{\rho}_n(i)=[g_n(i-1,i),g_n(i-2,i),\ldots,g_n(1,i)]^\texttt{T}$,  with $g_n(j,i)$ denoting the $j$th element of $\mathbf{g}_n(i)$.
   
Furthermore,   the associated error satisfies
   \begin{equation}\label{errotttrest}
 \varepsilon_{1n}(i)= \varepsilon_{1n}(i-1)\left[1-\lambda_n(i) (\lambda_n(i))^{(3-2n,1)}\right]
\end{equation}

\end{theorem}
\begin{proof}
Let $\mathcal{C}$ be the closed linear subspace associated with  $\mathbf{x}(1)-\hat{\mathbf{x}}_{n}(1,i)$. Denote by $ \boldsymbol{\phi}$ and $ \boldsymbol{\varphi}$  the projections of $\hat{\mathbf{x}}_{n}(i+1,i)$ onto $\mathcal{D}_2(i)$ and  $\mathcal{C}$, respectively\footnote{Subspace $\mathcal{C}$,  as well as the vectors $\boldsymbol{\phi}$ and $ \boldsymbol{\varphi}$, depend on $n$ and $i$; however,  for simplicity,  these parameters have been omitted.}.  Then,  in accordance with Theorem \ref{Teoremini}, it follows that
\begin{equation}\label{estunaeta}
\hat{\mathbf{x}}_{n}(i+1,i)=\boldsymbol{\phi}+\boldsymbol{\varphi}=\boldsymbol{\phi}+\lambda_n(i)\left(\mathbf{x}(1)-\hat{\mathbf{x}}_{n}(1,i)\right)
\end{equation}
with
\begin{equation*}
\lambda_n(i)=<\mathbf{x}(1)-\hat{\mathbf{x}}_{n}(1,i), \mathbf{x}(i+1)>_n \varepsilon_{2n}^{-1}(i)
\end{equation*}

By applying Proposition \ref{proppest} and using the Hermitian property of $\mathbf{T}^{(n)}$, equation  \eqref{coeficn} is derived.

Next, from \eqref{estunaeta} and Proposition \ref{proppest}, it follows that 
\begin{equation*}
\hat{\mathbf{x}}_{n}(i+1,i)=\lambda_n(i)\mathbf{x}(1)+\sum_{j=1}^{i-1}\left[f_n(j,i-1)+\lambda_n(i)g_n(i-j,i-1)
\right]\mathbf{x}(i+1-j)
\end{equation*}
where $f_n(j,i-1)$ is the $j$th element of $\mathbf{f}_n(i-1)$,   leading directly to equation \eqref{restodecoef}.

Finally, \eqref{errotttrest} follows from
\begin{multline*}
\varepsilon_{1n}(i)=\|\mathbf{x}(i+1)-\hat{\mathbf{x}}_{n}(i+1,i)\|_n^2=\|\mathbf{x}(i+1)-\boldsymbol{\phi}-\boldsymbol{\varphi}\|_n^2\\
=\|\mathbf{x}(i+1)-\boldsymbol{\phi}\|_n^2+\|\boldsymbol{\varphi}\|_n^2-<\mathbf{x}(i+1)-\boldsymbol{\phi},\boldsymbol{\varphi}>_n
\\  -<\boldsymbol{\varphi},\mathbf{x}(i+1)-\boldsymbol{\phi}>_n
\end{multline*}

\end{proof}
\section{Examples}\label{Examples}
\subsection{Theoretical Cases}

\subsubsection{Case 1}

%
This section illustrates the advantages of using an eigendecomposition based on a matrix different from the one proposed in \cite{Pei} for the tessarine case.  Moreover,  the convergence behavior of the power method is analyzed in a unique scenario.  Specifically,  consider the matrix

$$\mathbf{X}=\left[
  \begin{array}{ccc}
    x_1 & x_2 & x_3 \\
    x_2 & x_1 & x_4 \\
    x_3 & x_4 & x_1 \\
  \end{array}
\right]$$
with entries: $x_1=17+1.5\mathrm{i}+1.2\mathrm{j}+0.5\mathrm{k}$, $x_2=-0.1+0.2\mathrm{i}+0.03\mathrm{j}-0.04\mathrm{k}$,
$x_3=0.05+0.07\mathrm{j}-0.1\mathrm{k}$, and $x_4=0.05+0.1\mathrm{j}-0.2\mathrm{k}$.

Note that $\mathbf{X}$ is not a  $2$-Hermitian matrix. 
Consequently, expressing  $\mathbf{X}$
in terms of its eigenvalues and eigenvectors requires the use of equation \eqref{eigdes}, which involves computing the inverse of a matrix.  However, since $\mathbf{X}$ is $1$-Hermitian, it can be represented using equation \eqref{posher},  thereby avoiding this additional computational burden.

The approximations based on the eigenvalues and eigenvectors are analyzed in two settings.

Let $\lambda_n$ and $\mathbf{u}_n $ denote the eigenvalues and eigenvectors of $\mathbf{X}$ for some general parameters $\alpha, \beta>0$.  Similarly,  let $\varphi_n$ and $\mathbf{v}_n $ denote the eigenvalues and eigenvectors of $\mathbf{X}$ for the tessarine case ($\alpha=-1$, $ \beta=1$).   For $n=1,2$,  consider
$\mathbf{B}_n=\sum_{j=1}^n\lambda_j\mathbf{u}_j\mathbf{u}_j^{\texttt{H}_{1}^{1}}$
and $\mathbf{C}_n=\sum_{j=1}^n\varphi_j\mathbf{v}_j\mathbf{m}_j^\texttt{T}$, with $\mathbf{m}_j^\texttt{T}$  the $j$th  row of the inverse matrix that appears in equation \eqref{eigdes}.

Consider the extension of the real Frobenius norm for matrices in the  $(\alpha\beta)$-tessarine domain,  given by
$\|\mathbf{X}\|_F=\|\mathbf{X}\|_2$ for $\alpha=-1$, $\beta=1$.

Define $\pi_1(\alpha,\beta,n)=\|\mathbf{X}-\mathbf{C}_n\|_F-\|\mathbf{X}-\mathbf{B}_n\|_F$ and $\pi_2(\alpha,\beta,n)=\|\mathbf{X}-\mathbf{C}_n\|_1-\|\mathbf{X}-\mathbf{B}_n\|_1$.

\begin{figure}[h!]
\centering
\hskip-0.5cm \begin{minipage}[t]{7cm}
\hskip-0.5cm \includegraphics[width=8cm,height=4cm,draft=false]{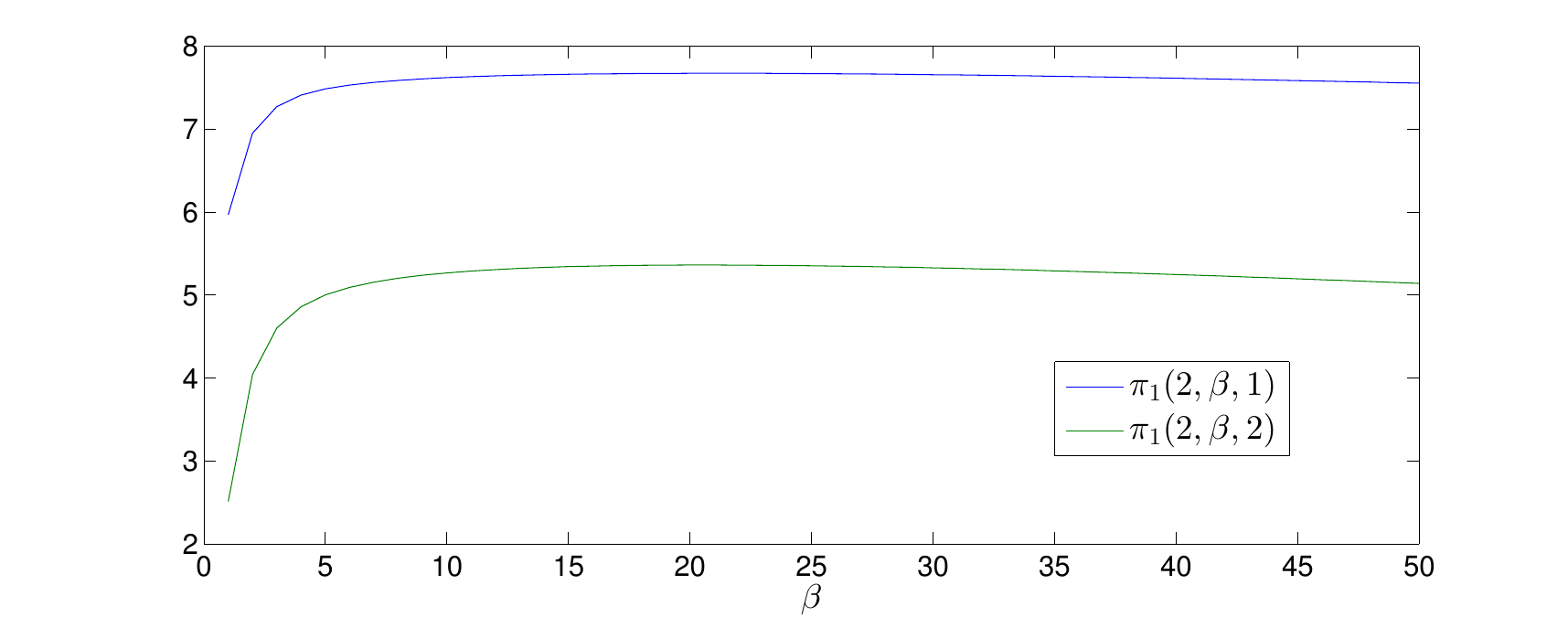}

\hskip-0.5cm \includegraphics[width=8cm,height=4cm,draft=false]{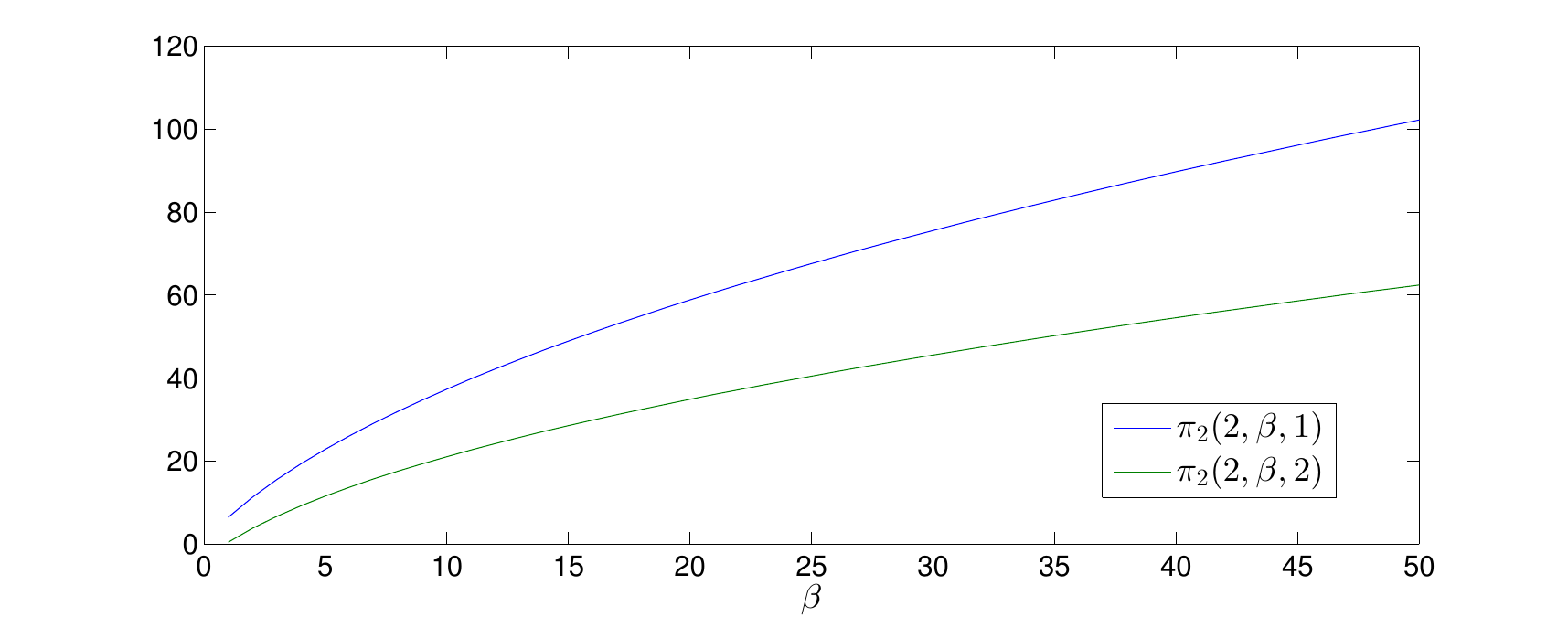}
\end{minipage}
\hskip-0.25cm \begin{minipage}[t]{7cm}
 \includegraphics[width=8cm,height=4cm,draft=false]{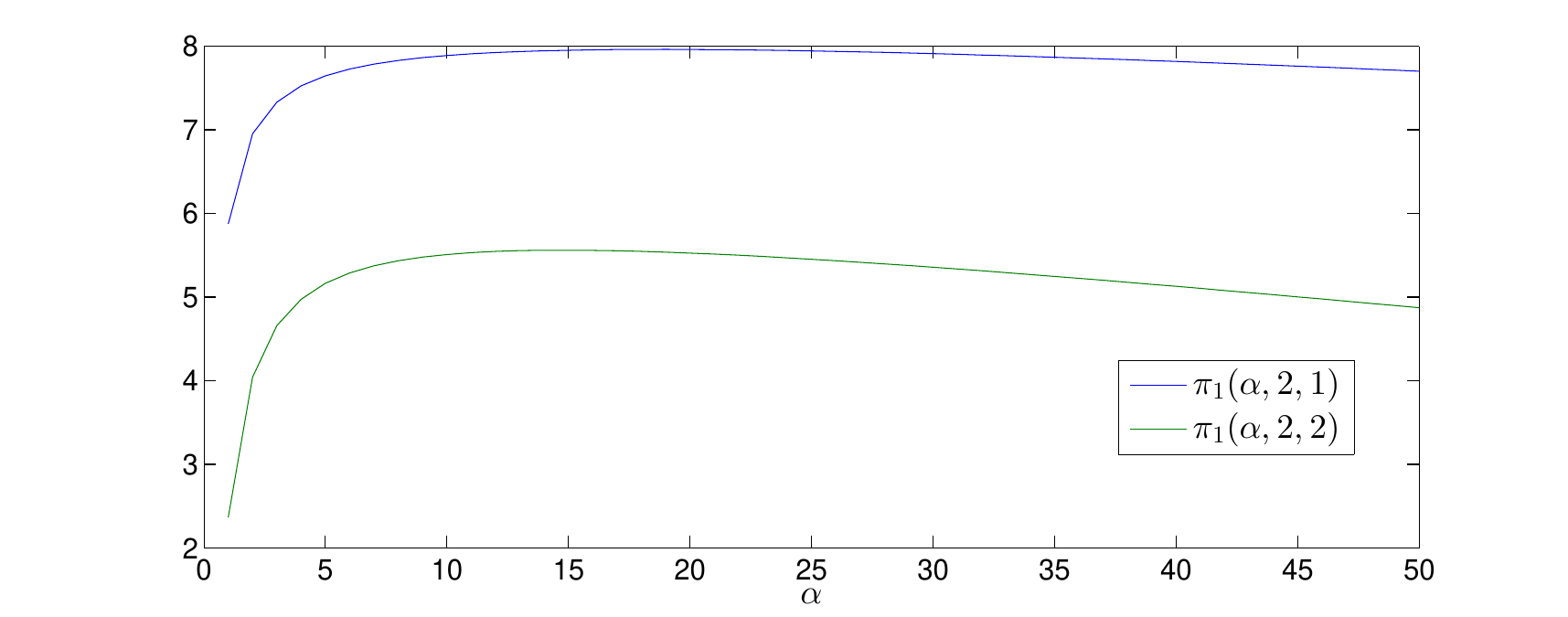} 
 \includegraphics[width=8cm,height=4cm,draft=false]{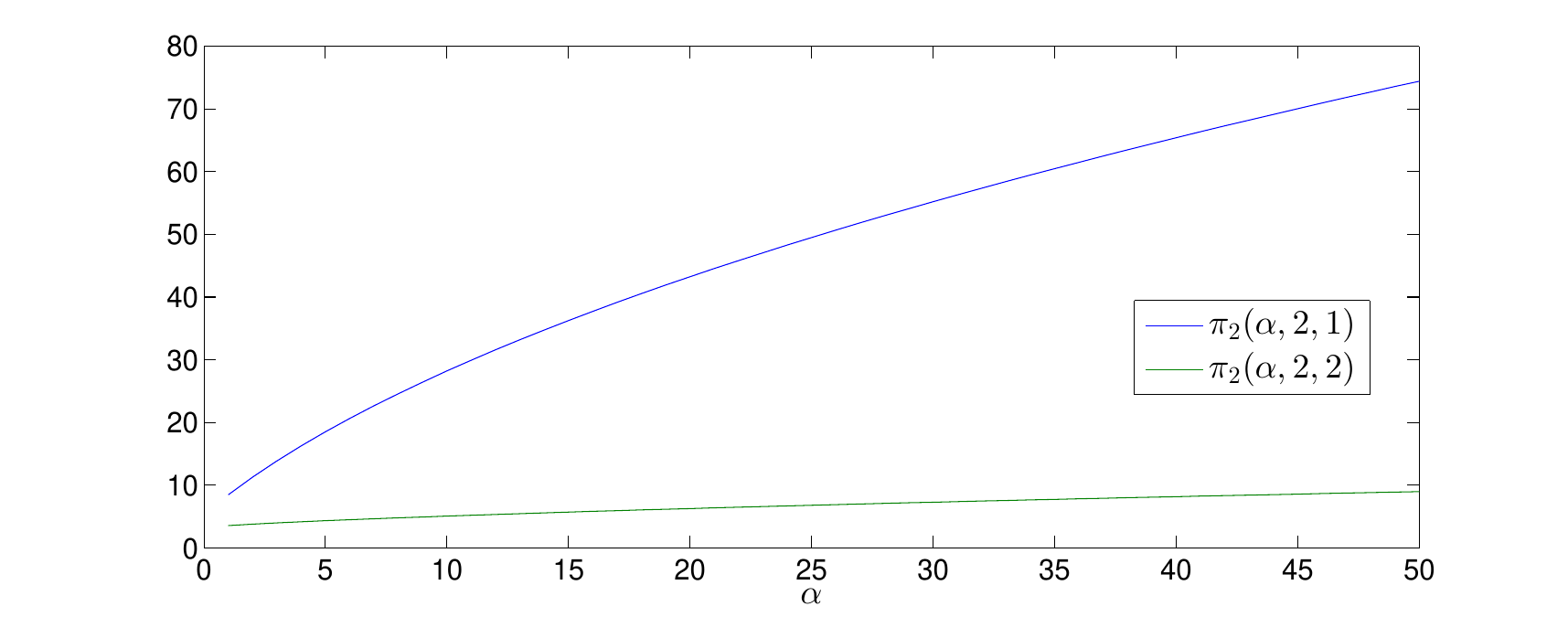}
\end{minipage} \hfill
\caption{$\pi_\nu(\alpha,\beta,n)$, $\nu=1,2$ and $n=1,2$.}\label{fig1}
\end{figure}

Figure \ref{fig1} analyzes two cases:
\begin{enumerate}
 \item $\alpha=2$ y $\beta=1,\ldots, 50$.
 \item $\beta=2$ y  $\alpha=1,\ldots, 50$.
\end{enumerate}

In all these scenarios, it can be observed that the approximation provided by $\mathbf{B}_n$ consistently outperforms that of $\mathbf{C}_n$, which is derived using the tessarine domain. While Corollary \ref{corolid} guarantees that $\pi_2(\alpha, \beta, n) \geq 0$, for all $\alpha, \beta$ and $n = 1,2$, it is noteworthy in this example that  $\pi_1(\alpha, \beta, n) \geq 0$ is also empirically satisfied for all $\alpha, \beta$ and $n = 1,2$.  The fact that $\pi_1(\alpha, \beta, n) \geq 0$ reflects that matrix $\mathbf{B}_n$ approximates matrix $\mathbf{X}$ better than $\mathbf{C}_n$, even though the $\|\cdot\|_F$ norm is the most appropriate in the tessarine domain. This example illustrates that, when the matrix $\mathbf{X}$ is not 2-Hermitian, the result on the best rank-$k$ approximation (established in the version of Corollary \ref{corolid} for $\alpha<0$) does not hold.  Furthermore, it can be seen that metrics $\pi_1$ and $\pi_2$ both exhibit similar behavior in the two cases examined, although a greater discrepancy is observed in $\pi_2$ than in $\pi_1$ as the values of $\alpha$ and $\beta$ increase.

Now,  a study of the power method is conducted.  First, case $\alpha = 14$ and $\beta = 2$ is analyzed, for which  the dominant eigenvalue is  $\lambda_1= 18.3227 +1.5082\mathrm{i}+1.1252\mathrm{j}+0.4873\mathrm{k}$. Second,  case $\alpha = 2$ and $\beta = 14$ is analyzed, where the dominant eigenvalue is $\lambda_1=18.1725 +1.2463\mathrm{i}+1.2045\mathrm{j}+0.4752\mathrm{k}$.

For this analysis, four $(\alpha\beta)$-vectors $\mathbf{y}_i$, $i = 1, \ldots, 4$, were randomly generated and  used as initial seeds  for the application of Theorem \ref{powerla},  in accordance with  expression
$\displaystyle \mathbf{x}_{0i}=\frac{\mathbf{y}_i}{|\mathbf{y}_i|_1}$, $i=1,\dots,4$.  The corresponding vectors $\mathbf{y}_1$ through $\mathbf{y}_4$ are:

\small{
\begin{equation*}
\mathbf{y}_1 = \left[
\begin{array}{c}
0.1352 - 0.9415\mathrm{i} - 0.532\mathrm{j} - 0.4838\mathrm{k} \\
0.5152 - 0.1623\mathrm{i} + 1.6821\mathrm{j} - 0.712\mathrm{k} \\
0.2614 - 0.1461\mathrm{i} - 0.8757\mathrm{j} - 1.1742\mathrm{k}
\end{array}
\right]
\end{equation*}
\begin{equation*}
\mathbf{y}_2 = \left[
\begin{array}{c}
0.3252-1.7115\mathrm{i}+0.3192\mathrm{j}-0.0301\mathrm{k} \\
-0.7549-0.1022\mathrm{i}+0.3129\mathrm{j}-0.1649\mathrm{k}\\
 1.3703 -0.2414\mathrm{i}-0.8649\mathrm{j}+0.6277\mathrm{k}
\end{array}
\right]
\end{equation*}
\begin{equation*}
\mathbf{y}_3 = \left[
\begin{array}{c}
 1.0933+0.0774\mathrm{i}-0.0068\mathrm{j}+0.3714\mathrm{k} \\
 1.1093-1.2141\mathrm{i}+1.5326\mathrm{j}-0.2256\mathrm{k} \\
-0.8637-1.1135\mathrm{i}-0.7697\mathrm{j}+1.1174\mathrm{k}
\end{array}
\right]
\end{equation*}
\begin{equation*}
\mathbf{y}_4 = \left[
\begin{array}{c}
0.7254 -0.205\mathrm{i} +1.409\mathrm{j} -1.2075\mathrm{k}\\
 -0.0631 -0.1241\mathrm{i}+1.4172\mathrm{j}+0.7172\mathrm{k}\\
 0.7147 +1.4897\mathrm{i}+0.6715\mathrm{j}+1.6302\mathrm{k}
\end{array}
\right]
\end{equation*}
}

Table \ref{tabla1PM} and Table \ref{tabla2PM} report the approximated dominant eigenvalues after $30$ and $100$ iterations for each initialization. The results confirm convergence to the expected dominant eigenvalues, with increasingly accurate approximations over successive iterations.

\begin{table}[ht]
\centering
\footnotesize{
\begin{tabular}{|c|c|c|}
  \hline
  \multicolumn{2}{|c|}{\(\alpha = 14, \quad \beta = 2\)} & {\color{red}\(\lambda_1 = 18.3227 + 1.5082\mathrm{i} + 1.1252\mathrm{j} + 0.4873\mathrm{k}\)} \\
  \hline
 & 30 & 100 \\
  \hline
  \(\mathbf{x}_{01}\) & \(17.9785 + 1.4246\mathrm{i} + 0.9044\mathrm{j} + 0.4227\mathrm{k}\) & \(18.3218 + 1.508\mathrm{i} + 1.1245\mathrm{j} + 0.4871\mathrm{k}\) \\
  \hline
    \(\mathbf{x}_{02}\) & \(18.1911+1.5392\mathrm{i}+1.033\mathrm{j}+0.5094\mathrm{k}\) & \(18.3227 +1.5082\mathrm{i}    +1.1252\mathrm{j}+0.4873\mathrm{k}\) \\
  \hline
      \(\mathbf{x}_{03}\) & \(17.9255 +1.4103\mathrm{i}+0.8565\mathrm{j}+0.415\mathrm{k}\) & \(18.3109+1.5051\mathrm{i}    +1.1168\mathrm{j}+0.485\mathrm{k}\) \\
  \hline
        \(\mathbf{x}_{04}\) & \(18.2776 +1.5099\mathrm{i} +1.1316\mathrm{j}+0.4801\mathrm{k}\) & \(18.3227 +1.5082\mathrm{i}+1.1252\mathrm{j}+0.4873\mathrm{k}\) \\
  \hline
\end{tabular}
}
\caption{Case \(\alpha = 14\) and \(\beta = 2\).}
\label{tabla1PM}
\end{table}

\begin{table}[ht]
\centering
\footnotesize{
\begin{tabular}{|c|c|c|}
  \hline
  \multicolumn{2}{|c|}{\(\alpha = 2, \quad \beta = 14\)} & {\color{red}\(\lambda_1=18.1725 +1.2463\mathrm{i}+1.2045\mathrm{j}+0.4752\mathrm{k}\)} \\
  \hline
  & 30 & 100 \\
  \hline
  \(\mathbf{x}_{01}\) & \( 18.1248+1.221\mathrm{i}+1.195\mathrm{j}+0.4664\mathrm{k}\) & \(18.1717+1.2457\mathrm{i}    +1.2043\mathrm{j}+0.475\mathrm{k}\) \\
  \hline
    \(\mathbf{x}_{02}\) & \( 18.1521+1.2437\mathrm{i}+1.2014\mathrm{j}+0.4761\mathrm{k}\) & \(18.1723+1.2461\mathrm{i}     +1.2045\mathrm{j}+0.4752\mathrm{k} \) \\
  \hline
   \(\mathbf{x}_{03}\) & \( 17.7518+1.5033\mathrm{i}+1.0922\mathrm{j}+0.5439\mathrm{k}\) & \(18.1722 +1.2461\mathrm{i}   +1.2045\mathrm{j}+0.4752\mathrm{k}
 \) \\
  \hline
   \(\mathbf{x}_{04}\) & \(17.6044+1.4434\mathrm{i}+1.0667\mathrm{j}+0.5184\mathrm{k}\) & \( 18.1721    +1.246\mathrm{i}+1.2044\mathrm{j}+0.4751\mathrm{k}\) \\
  \hline
\end{tabular}
}
\caption{Case \(\alpha = 2\) and \(\beta = 14\).}
\label{tabla2PM}
\end{table}


\subsubsection{Case 2}

This section presents a comparative analysis of the execution times of the algorithms described in Section \ref{sectiolineal} for computing projections in least squares problems. Specifically,  the computational advantages of the sequential algorithm introduced in Theorem \ref{secuencialalgort}, particularly when applied to Toeplitz $(\alpha\beta)$-tessarine matrices, are highlighted in contrast to more conventional approaches.

For this purpose, considering that $\alpha = -2$ and $\beta = 3$,  a square matrix $\mathbf{X}$ of order $500$ was generated, such that the associated matrix $\mathbf{T}^{(2)}$, as defined in equation \eqref{matrixT}, is a $2$-Hermitian Toeplitz matrix.  This matrix was generated using the built-in MATLAB function specifically designed to create this type of Toeplitz matrix.  Then, fixed-point updating projection $\mathbf{x}_2(1,i)$ was computed for $i=1,\dots, 500$ using the following three methods: the inverse method, as described in equation \eqref{parminv};  the pseudoinverse method, based on equation \eqref{leastpseuuu}; and the sequential algorithm, as presented in Theorem \ref{secuencialalgort}.

\begin{figure}[h!]
\centering{
\includegraphics[width=10.5cm,height=8.2cm]{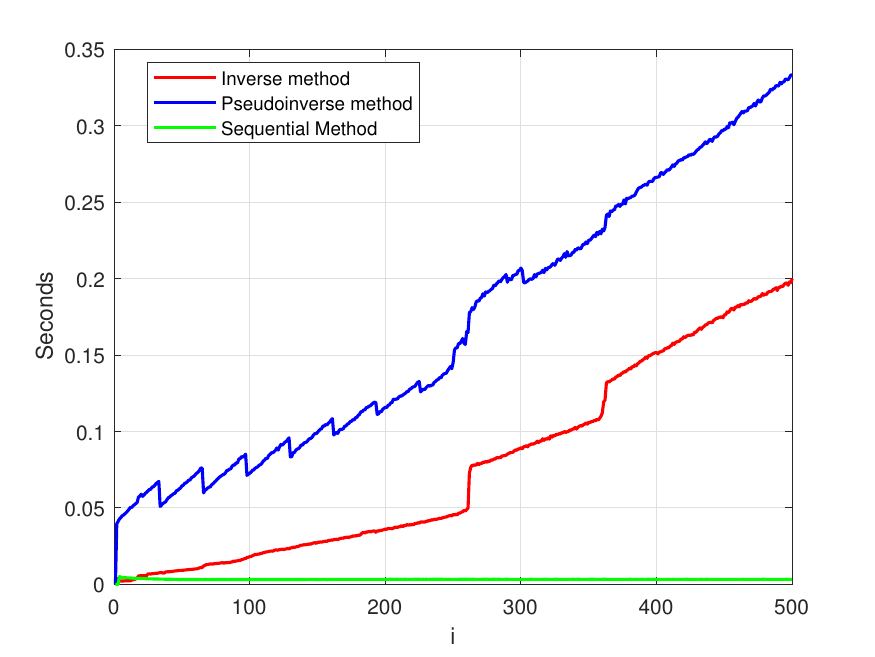}}\caption{Average computation time (in seconds) for the estimation  of $\mathbf{x}_2(1,i)$  using the inverse method (red line), the pseudo inverse method (blue line), and the sequential algorithm (green line).\label{tiempocostp}
}\end{figure}

Figure \ref{tiempocostp} displays the average computation times (in seconds) obtained over  $200$ simulations for each of the three methods.  All computations were performed on a standard laptop equipped with  a 13th Gen Intel(R) Core(TM) i9-13980HX 2.20 GHz processor and 32 GB of RAM.

The results reveal that the sequential algorithm has a clear computational advantage when the structural conditions required for its application are satisfied. Moreover, the inverse and pseudoinverse methods exhibit exponential growth in average computation time as $i$ increases, whereas the sequential algorithm maintains a linear growth trend. This behavior underscores the superior efficiency and scalability of the sequential method for Toeplitz matrices.

It is worth emphasizing that the specific values of parameters $\alpha$ and $\beta$ do not significantly influence the computational behavior observed in the study. Consistent results have been obtained across different $(\alpha\beta)$-tessarine algebras, suggesting that the findings presented are robust and broadly applicable within the algebraic framework considered.

\subsection{Real Applications}
\subsubsection{Application of SVD to the  Analysis of the Watermarked Host Image}

In this example, image reconstruction based on SVD is examined, with a comparative analysis between $(\alpha\beta)$-tessarine and standard tessarine algebras.

To this end, the process begins with two original images denoted as $\mathbf{A}$ and $\mathbf{B}$. Image $\mathbf{A}$ is considered as the host image, while image $\mathbf{B}$ is the watermark. A collage image, referred to as the watermarked host image and denoted by $\mathbf{A}_{\mathbf{B}}=\mathbf{A}+\mu\mathbf{B}$, is then constructed for $\mu=0.04,0.06,0.08,0.1$.

Subsequently, image quality is assessed by comparing the original images with their approximated counterparts obtained from the watermarked host image using the SVD. For this purpose, the peak signal-to-noise ratio (PSNR) metric is employed, defined as 
\begin{equation*}
\text{PSNR}=10\log_{10}\left(\frac{255 ^2}{\frac{MSE_R+MSE_G+MSE_B}{3}}\right)
\end{equation*}
where $MSE_\nu$, for $\nu=R,G,B$, denotes the mean squared error between the pixel values of color channel $\nu$ in the original image and those in the approximated image.

\begin{figure}
\centering
\includegraphics[width=\textwidth]{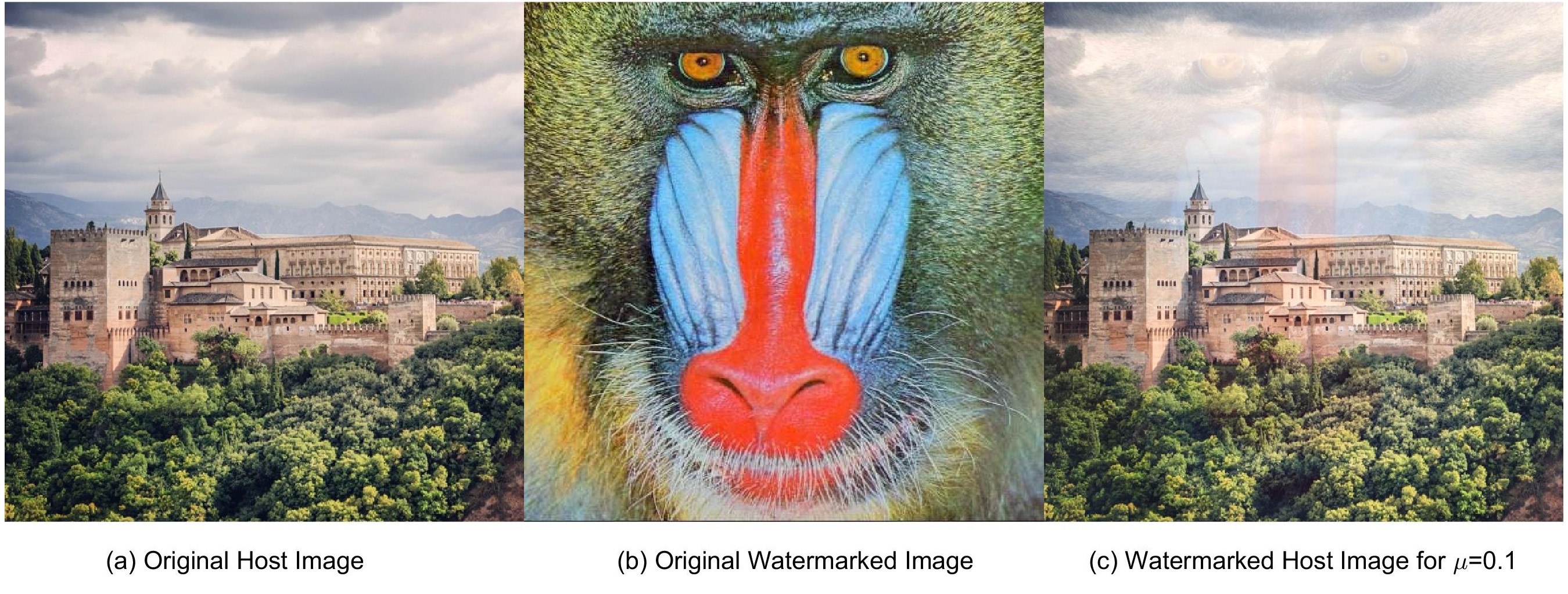}
\caption{(a) Image 1: ``{\em The Alhambra Palace}''; (b) Image 2: ``{\em Mandrill}''; (c) Watermarked Host Image for $\mu=0.1$. }
\label{watermarked Host Image 1}
\end{figure}

\begin{figure}
\centering
\includegraphics[width=\textwidth]{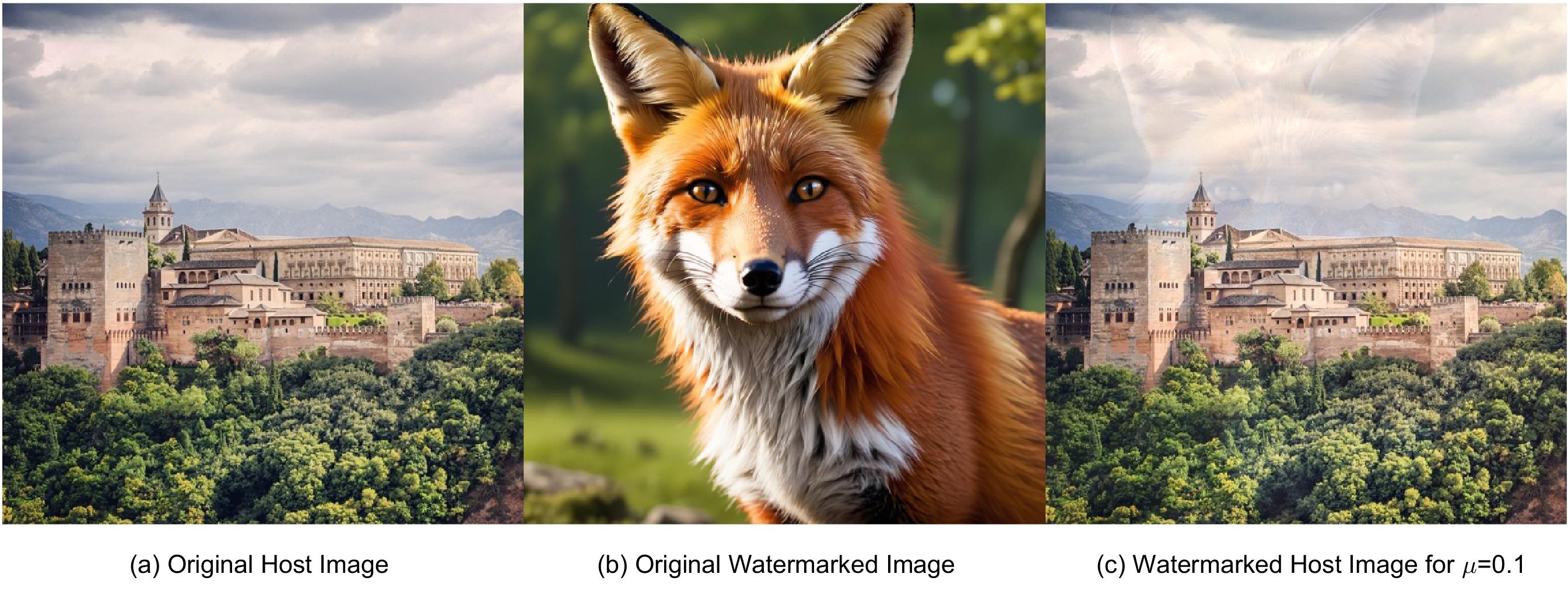}
\caption{(a) Image 1: ``{\em The Alhambra Palace}''; (b) Image 3: ``{\em Fox}''; (c) Watermarked Host Image for $\mu=0.1$.  }
\label{watermarked Host Image 2}
\end{figure}

This analysis is conducted under two distinct configurations. These configurations are illustrated in Figures \ref{watermarked Host Image 1} and \ref{watermarked Host Image 2} , which display the original images, $\mathbf{A}$ and $\mathbf{B}$, as well as the resulting watermarked host image constructed from them, $\mathbf{A}_B$, with $\mu = 0.1$.

The images used in this study are either openly available on Flickr or are commonly employed in image processing experiments:
\begin{itemize}
\item  Image 1: ``{\em The Alhambra Palace}''. Image under license CC BY 2.0, source: \url{https://www.flickr.com/photos/ichstyle/19212468916/in/dateposted/}.
\item Image 2: ``{\em Mandrill}''. Source: USC-SIPI Image Database,  available at: \url{http://sipi.usc.edu/database/database.php?volume=misc\&image=10\#top}. 
\item Image 3: ``{\em Fox}''. Image under license CC PDM 1.0, source: \url{https://www.flickr.com/photos/ioachimphotos/53485096097/in/photolist-2puotuY-2pui6Fg-2puotve-2puothZ-2punK1x}.
\end{itemize}

Moreover, for each of these configurations, PSNR values for the host image, $\mathbf{A}$, have been computed using both the tessarine SVD (in black) and the $(\alpha\beta)$-tessarine SVD with parameters $\alpha=3$, $\beta=1$ (in red), considering 50, 150, and 250 singular values. These results are reported in Tables \ref{PSNR host image 1} and \ref{PSNR host image 2} for the scenarios depicted in Figures 5 and 6, respectively. A similar analysis has been carried out for the watermarked image, $\mathbf{B}$, this time considering 150, 250, and 350 singular values. The corresponding results are presented in Tables \ref{PSNR watermarked image 1} and \ref{PSNR watermarked image 2} for each setting.

The results obtained show that, on the one hand, the PSNR values increase as the number of singular values used for reconstruction grows, and on the other hand, that they decrease as the value of $\mu$ increases. However, the most significant conclusion to be drawn is that the PSNR values obtained using the $(\alpha\beta)$-tessarine SVD (shown in red), are consistently higher than those obtained with the tessarine SVD (shown in black). This example therefore illustrates the superiority of the $(\alpha\beta)$-tessarine SVD approach over that of the tessarine SVD.

\begin{table}
\begin{center}
\begin{tabular}{|c|c|c|c|}
  \hline

   & 50 & 150 & 250 \\
    \hline
  0.04 & (25.3749, {\color{red}25.4443}) &(30.8842, {\color{red}30.9364})  & (32.9609, {\color{red}32.9720})\\
  0.06 & (24.5929, {\color{red}24.6495}) &(28.5519, {\color{red}28.5815})& (29.6300, {\color{red}29.6352}) \\
  0.08 & (23.6912, {\color{red}23.7360})& (26.5601, {\color{red}26.5781}) & (27.1990, {\color{red}27.2020})\\
  0.1 & (22.7530, {\color{red}22.7881 })& (24.8777, {\color{red}24.8897}) & (25.2916, {\color{red}25.2934})\\
  \hline
\end{tabular}
\caption{PSNR for the approximated  host image (Figure \ref{watermarked Host Image 1}).}
\label{PSNR host image 1}
\end{center}
\end{table}

\begin{table}
\begin{center}
\begin{tabular}{|c|c|c|c|}
  \hline

   & 50 & 150 & 250 \\
    \hline
  0.04 & (25.6028, {\color{red}25.6746}) &(31.7694, {\color{red}31.8351})  & (34.5117, {\color{red}34.5248})\\
  0.06 & (25.0311, {\color{red}25.0931}) &(29.7582, {\color{red}29.7994})& (31.2658, {\color{red}31.2726}) \\
  0.08 & (24.3381, {\color{red}24.3905})& (27.9419, {\color{red}27.9689}) & (28.8675, {\color{red}28.8713})\\
  0.1 & (23.5822, {\color{red}23.6272 })& (26.3593, {\color{red}26.3777}) & (26.9761, {\color{red}26.9785})\\
  \hline
\end{tabular}
\caption{PSNR for the approximated  host image (Figure \ref{watermarked Host Image 2}).}
\label{PSNR host image 2}
\end{center}
\end{table}

\begin{table}
\begin{center}
\begin{tabular}{|c|c|c|c|}
  \hline

   & 150 & 250 & 350 \\
    \hline
  0.04 & (6.5141, {\color{red}6.6323}) &(15.8340, {\color{red}15.9654})  & (29.7266, {\color{red}29.8142})\\
  0.06 & (9.9629, {\color{red}10.0752}) &(19.0611, {\color{red}19.1873})& (31.6817, {\color{red}31.7515}) \\
  0.08 & (12.3630, {\color{red}12.4688})& (21.1924, {\color{red}21.3079}) & (32.7978, {\color{red}32.8810})\\
  0.1 & (14.1798, {\color{red}14.2806})& (22.7218, {\color{red}22.8288}) & (33.5606, {\color{red}33.6416})\\
  \hline
\end{tabular}
\caption{PSNR for the approximated watermarked image (Figure \ref{watermarked Host Image 1}).}
\label{PSNR watermarked image 1}
\end{center}
\end{table}

\begin{table}
\begin{center}
\begin{tabular}{|c|c|c|c|}
  \hline

   & 150 & 250 & 350 \\
    \hline
  0.04 & (6.5557, {\color{red}6.6795}) &(16.0191, {\color{red}16.1329})  & (31.3846, {\color{red}31.4282})\\
  0.06 & (10.0510, {\color{red}10.1743}) &(19.4404, {\color{red}19.5664})& (34.2708, {\color{red}34.3240}) \\
  0.08 & (12.5135, {\color{red}12.6333})& (21.8082, {\color{red}21.9353}) & (36.0602 {\color{red}36.1223})\\
  0.1 & (14.4070, {\color{red}14.5287})& (23.5934, {\color{red}23.7190}) & (37.2845, {\color{red}37.3588})\\
  \hline
\end{tabular}
\caption{PSNR for the approximated watermarked image (Figure \ref{watermarked Host Image 2}).}
\label{PSNR watermarked image 2}
\end{center}
\end{table}

\subsubsection{Application of the Least Square Problem to Image Classification}

This section explores the application  of the $(\alpha\beta)$-tessarine least squares approach proposed  to the task of color face recognition,  formulated as a linear regression classification (LRC) problem.  
  
This is a benchmark example which has been frequently used in prior studies to evaluate  the performance of hypercomplex regression models.  Within this framework, the problem was initially addressed by \cite{Melegy},  who demonstrated the efficacy of tessarine representations in classification tasks,  highlighting their superiority over conventional quaternion-based methods.  The LRC approach was subsequently employed by \cite{JC} in the quaternion domain,  further reinforcing its importance as a comparative standard.  

Extending  these foundational contributions,  an enhanced solution based on the $(\alpha\beta)$-tessarine algebra is presented here,  which outperforms the one proposed in \cite{Melegy} based on tessarines. For this purpose,  the methodology introduced in \cite{Melegy} is established as a baseline for comparative analysis against the solutions presented  in the $(\alpha\beta)$-tessarine domain.  

Experimental validation is performed  on the FEI color face database \cite{Melegy,  JC, Thomaz}. This dataset comprises 2800 images of 200 individuals, captured against a uniform white background. The images exhibit variations in facial expressions and viewing angles. In this experiment, a subset of $L = 50$ subjects has been randomly selected. For each subject, $n = 7, 8, 9,$ or $10$ images have been used for training, while 3 different images have been reserved for testing. All the images have been down-sampled using a reduction factor of $0.1$, consistent with the preprocessing steps described in \cite{Melegy}.

The classification approach employed involves predicting the identity of each individual based on the regression output corresponding to their three test images, using the $n$ training images as input. The accuracy of this classification is measured as the percentage of correct identified individuals across all test samples.

Table \ref{accuracy} reports the recognition accuracy of both methods under comparison. As expected, the accuracy of both methods increases with the number of training images. Moreover, the $(\alpha\beta)$-tessarine configuration with  $\alpha=1.5$ and $\beta=1$ consistently outperforms the approach in \cite{Melegy}  (corresponding to the case $\alpha=-1$, $\beta=1$)  across all training sizes, indicating potential advantages of this algebraic framework in color image recognition tasks.

\begin{table}
\begin{center}
\begin{tabular}{|c|c|c|c|c|}
  \hline

 n  & 7 & 8 & 9& 10 \\
    \hline
 $ \alpha=-1$, $\beta=1$ & 90\% & 92\% & 92\%&$92.\wideparen{6}\%$\\
   $\alpha=1.5$, $\beta=1$ & {\color{red}$90.\wideparen{6}\%$} & {\color{red}$92.\wideparen{6}\%$} &  {\color{red}$92.\wideparen{6}\%$} &{\color{red}$93.\wideparen{3}\%$}  \\
  \hline
\end{tabular}
\caption{Accuracy rate.}\label{accuracy}
\end{center}
\end{table}

\section{Conclusions}\label{section Conclusions}

This paper  presents a comprehensive suite of algebraic tools specifically designed for $\alpha\beta$-tessarine matrices,  where  $\alpha \in \mathbb{R}-\{0\}, \beta > 0$.  These tools comprise methodologies for computing matrix inverses, square roots,  LU factorization with partial pivoting,  determinants,  eigenvalues and eigenvectors,  the power method,  SVD,  rank-$k$ approximations,  pseudoinverses,  and  solutions to  least squares problems. Notably,  this algebraic framework extends GSQ algebra \cite{Navarro_Segre} (case  $\beta = 1$)   and, consequently, elliptic quaternion algebra \cite{Kosal1}  (case $\alpha < 0$, $\beta = 1$).

The mathematical tools developed within this research are highly significant, providing a solid basis for prospective research, and they have a variety of applications across diverse scientific and engineering disciplines. These tools are not arbitrary; rather, they have been carefully created to meet the analytical demands of contemporary signal processing techniques. They form the essential theoretical foundation for the practical implementation of signal processing tasks,  enabling the development of efficient, robust, and scalable computational solutions. Therefore, their design responds directly to the structural and algorithmic needs of modern applications, especially in areas where multidimensional data and hypercomplex representations are involved.

Their particular relevance in signal processing is underscored by several key contributions.  In this context,  the degree of properness of a signal —based on some second-order statistics vanishing —  plays a critical role. This property has important implications for signal modeling and processing,  as it contributes to a significant reduction in computational complexity,  which is essential for the efficient handling of high-dimensional signals. The algebraic tools proposed are thus particularly well suited to exploit this structure, offering both theoretical rigor and practical efficiency.  

For instance,  the matrix inversion technique is essential for formulating and solving signal estimation and detection problems,  while preserving the structural advantages offered by properness  \cite{Navarro_beta_quaternion,  Bihan,  Wang,  Krupinski}.  When matrix inversion is not feasible,  the matrix pseudoinverse serves as a key alternative,  enabling least squares solutions in linear regression problems with ill-conditioned matrices  \cite{JC}.  It also plays a fundamental role in SVD, which underpins modern techniques in image compression and dimensionality reduction.  In particular, its applications in Principal Component Analysis  and Partial Least Squares make it indispensable for analyzing complex data structures and improving signal processing techniques \cite{JC}.

Computations of matrix square roots and determinants are also key elements in numerous signal processing applications.  In particular,  the determinant becomes especially relevant in detection problems, enabling generalized likelihood ratios, which are essential in statistical tests designed to discriminate between hypotheses involving proper and improper signals \cite{Bihan,  Krupinski}, to be evaluated.  

In addition,  the LU factorization with partial pivoting technique constitutes a robust method for solving  linear systems,   computing matrix inverse and determinants,  while also improving  numerical stability.   It is particularly effective for handling singular or ill-conditioned matrices in adaptive filtering, dynamic system modeling, and spectral analysis \cite{JC}.

In the context of spectral matrix decompositions, the power method  provides a powerful and efficient approach for extracting the dominant eigenvalue and  eigenvector in spectral matrix decompositions. This method can be widespread applied in signal processing, data compression, optimization, and stability analysis in dynamic systems \cite{Li, Cui}.  

Finally,  the inclusion of Toeplitz matrices enables stationary properties of signals to be captured, as a result of their autocovariance structure. This fact facilitates the analysis and modeling of stationary signals via efficient spectral decomposition. The structured representation of Toeplitz matrices is widely used in filtering,  spectral analysis, signal prediction, and data compression. \cite{Brockwell, Huang}.

The computational implementation of these methodologies is essential for practical use. 
For this purpose,  MATLAB-based functions have been  developed to support the numerical experiments presented in Section~\ref{Examples}.  These developments are the first step toward the creation of a dedicated MATLAB toolbox, designed to facilitate broader accessibility  and usability of the $\alpha\beta$-tessarine framework in both academic and applied research settings. Future efforts will focus on the optimization and expansion of  these computational tools,  and the explorations of their practical applications in various signal processing scenarios.

\section{CRediT Authorship Contribution Statement}

{\bf José Domingo Jiménez-López}: Writing – original draft, Methodology, Investigation, Conceptualization, Software. {\bf Jesús Navarro-Moreno}: Writing – original draft, Methodology, Investigation, Conceptualization, Software. {\bf Rosa María Fernández-Alcalá}: Writing – review \& editing, Formal analysis, Funding acquisition. {\bf Juan Carlos Ruiz-Molina}: Validation, Supervision, Formal analysis, Visualization.

\section{Conflict of Interest Declaration}

The authors declare that they have no known conflicting financial interests or personal relationships that could appearto influence the work presented in this paper.

\section{Acknowledgments}

This work is part of the I+D+i project PID2021-124486NB-I00, funded by
MICIU/AEI/10.13039/501100011033/ and ERDF/EU.

\appendix
\section{Methods for $(\alpha)$-Complex Matrix Transformations}\label{Ap1}

This appendix presents computational methods for performing transformations on 
$(\alpha)$-complex matrices,  including the inverse, square root, and LU factorization with partial pivoting.   Additionally,  approaches for determining the eigendecomposition and the SVD of these matrices are presented.


Throughout this section,  the matrices  $\mathbf{X}=\mathbf{A}+\mathbf{B}\mathrm{i} \in \mathbb{C}_{\alpha}^{p\times q}$,  where  $\mathbf{A}, \mathbf{B}\in \mathbb{R}^{p\times q}$,  and  $\dot{\mathbf{X}}=\mathbf{A}+\dot{\mathbf{B}}\epsilon \in \mathbb{C}^{p\times q}$,   with $\dot{\mathbf{B}}=\sqrt{|\alpha|}\mathbf{B}$ and $\epsilon^2=-1$,  are considered. 
  \subsection{Inverse  of  an  $(\alpha)$-Complex Matrix}
  
  The following Lemma  provides a way for computing the inverse of an 
$(\alpha)$-complex matrix.
  

\begin{lemma}\label{lema1}
Let $\mathbf{X}\in \mathbb{C}_{\alpha}^{p\times p}$.   The following properties hold:
\begin{description}
\item[ \normalfont{Case }$\alpha<0$:] $\mathbf{X}$ is invertible  if, and only if, $\dot{\mathbf{X}}$ is invertible.  Furthermore,  if the inverse of $\dot{\mathbf{X}}$ exists and takes the form $\dot{\mathbf{X}}^{-1}=\mathbf{C}+\dot{\mathbf{D}}\epsilon$,   then the inverse of $\mathbf{X}$ is given by $\mathbf{X}^{-1}=\mathbf{C}+\mathbf{D}\mathrm{i}$,  with $\mathbf{D}=\dot{\mathbf{D}}/\sqrt{|\alpha|}$. 
    \item[ \normalfont{Case }$\alpha>0$:] Consider the real matrices $\mathbf{E}=(\mathbf{A}+\dot{\mathbf{B}})^{-1}$ and $\mathbf{F}=(\mathbf{A}-\dot{\mathbf{B}})^{-1}$. Then,   $\mathbf{X}$  is invertible if, and only if, both $\mathbf{E}^{-1}$ and $\mathbf{F}^{-1}$ exist.  Furthermore,  if $\mathbf{X}$ is invertible,  its inverse is given by   $\mathbf{X}^{-1}=\mathbf{C}+\mathbf{D}\mathrm{i}$,  with $\mathbf{C}=\frac{\mathbf{E}+\mathbf{F}}{2}$ and $\mathbf{D}=\frac{\mathbf{E}-\mathbf{C}}{\sqrt{\alpha}}$.
\end{description}
\end{lemma}

\subsection{Square Root of an  $(\alpha)$-Complex Matrix}

  The result below proposes a method for computing the square root  of an 
$(\alpha)$-complex matrix.


\begin{lemma}\label{lema221}
Let $\mathbf{X}\in \mathbb{C}_{\alpha}^{p\times p}$.  The following properties hold:
\begin{description}
\item[ \normalfont{Case }$\alpha<0$:] The square root of $\mathbf{X}$  exists if, and only if,  the square root of $\dot{\mathbf{X}}$  exists.  Furthermore,  if the square root of $\dot{\mathbf{X}}$  exists and is given by $\dot{\mathbf{X}}^{1/2}=\mathbf{C}+\dot{\mathbf{D}}\epsilon$,   then  the square root of $\mathbf{X}$ is given by $\mathbf{X}^{1/2}=\mathbf{C}+\mathbf{D}\mathrm{i}$,  with $\mathbf{D}=\dot{\mathbf{D}}/\sqrt{|\alpha|}$.
    \item[ \normalfont{Case }$\alpha>0$:] Consider  the real matrices $\mathbf{E}=(\mathbf{A}+\dot{\mathbf{B}})^{1/2}$ and $\mathbf{F}=(\mathbf{A}-\dot{\mathbf{B}})^{1/2}$. Then,  the square root of $\mathbf{X}$  exists if, and only if, both $\mathbf{E}^{1/2}$ and $\mathbf{F}^{1/2}$  exist.   Furthermore,  in this case, the square root of $\mathbf{X}$  is given by  $\mathbf{X}^{1/2}=\mathbf{C}+\mathbf{D}\mathrm{i}$, with $\mathbf{C}=\frac{\mathbf{E}+\mathbf{F}}{2}$ and $\mathbf{D}=\frac{\mathbf{E}-\mathbf{C}}{\sqrt{\alpha}}$.
\end{description}
\end{lemma}

\subsection{LU Factorization with Partial Pivoting of an $(\alpha)$-Complex Matrix}
The following lemma outlines the LU factorization  with partial pivoting  for $(\alpha)$-complex matrices. 
\begin{lemma}\label{LUlema}
Consider $\mathbf{X}\in \mathbb{C}_{\alpha}^{p\times p}$.
  The following properties hold:
\begin{description}
\item[ \normalfont{Case }$\alpha<0$:] Let   $(\dot{\mathbf{P}},\dot{\mathbf{L}}, \dot{\mathbf{U}})$    be  the matrices of the LU factorization with partial pivoting  of the complex matrix $\dot{\mathbf{X}}$,  where $ \dot{\mathbf{P}}$ is a real matrix, and $\dot{\mathbf{L}}, \dot{\mathbf{U}}$ are complex matrices that take the form
    \begin{equation*}\label{matneds}
    \dot{\mathbf{L}}={\mathbf{L}}_1+{\mathbf{L}}_2\epsilon,\quad
    \dot{\mathbf{U}}={\mathbf{U}}_1+{\mathbf{U}}_2\epsilon
     \end{equation*}
  Then, the matrices  $(\mathbf{P},\mathbf{L},\mathbf{U})$ of the LU factorization with partial pivoting  of the matrix $\mathbf{X}$ are given by
        \begin{equation*}
     \displaystyle{\mathbf{P}=\dot{\mathbf{P}}}, \quad \displaystyle{\mathbf{L}={\mathbf{L}}_1+\frac{{\mathbf{L}}_2} {\sqrt{|\alpha|}}\mathrm{i}},\quad 
    \displaystyle{\mathbf{U}={\mathbf{U}}_1+\frac{{\mathbf{U}}_2}{\sqrt{|\alpha|}}\mathrm{i}}
     \end{equation*}
         \item[ \normalfont{Case }$\alpha>0$:] Let
    \begin{equation}\label{matrds}
    (\mathbf{P}_s,\mathbf{L}_s,\mathbf{U}_s), \quad(\mathbf{P}_d,\mathbf{L}_d,\mathbf{U}_d)
       \end{equation}
       be the matrices of the LU factorization with partial pivoting  of the real matrices  $\mathbf{A}+\dot{\mathbf{B}}$ and $\mathbf{A}-\dot{\mathbf{B}}$, respectively. Then, the matrices of the LU factorization with partial pivoting  of matrix $\mathbf{X}$ take the form
        \begin{equation*}
        \begin{split}
        \mathbf{P}&=\frac{\mathbf{P}_s+\mathbf{P}_d}{2}+\frac{\mathbf{P}_s-\mathbf{P}_d}{2\sqrt{\alpha}}\mathrm{i}\\
         \mathbf{L}&=\frac{\mathbf{L}_s+\mathbf{L}_d}{2}+\frac{\mathbf{L}_s-\mathbf{L}_d}{2\sqrt{\alpha}}\mathrm{i}\\
          \mathbf{U}&=\frac{\mathbf{U}_s+\mathbf{U}_d}{2}+\frac{\mathbf{U}_s-\mathbf{U}_d}{2\sqrt{\alpha}}\mathrm{i}
        \end{split}
        \end{equation*}
\end{description}
\end{lemma}

\subsection{Eigendecomposition of an $(\alpha)$-Complex Matrix}

The following lemma provides explicit formulas to calculate the eigenvalues and eigenvectors of an $(\alpha)$-complex matrix.  

\begin{lemma}\label{lema2}
Consider that $\mathbf{X}\in \mathbb{C}_{\alpha}^{p\times p}$. 
The following properties hold:
\begin{description}
\item[ \normalfont{Case }$\alpha<0$:] Let   $\dot{\lambda}={\lambda}_1+{\lambda}_2\epsilon$ and $\dot{\mathbf{u}}=\mathbf{u}_1+\mathbf{u}_2\epsilon$ be  the eigenvalues and eigenvectors of complex matrix $\dot{\mathbf{X}}$. Then, the eigenvalues and eigenvectors of $\mathbf{X}$ are given by $\lambda={\lambda}_1+\displaystyle \frac{{\lambda}_2}{\sqrt{|\alpha|}}\mathrm{i}$ and
    $\mathbf{u}=\mathbf{u}_1+\displaystyle \frac{\mathbf{u}_2}{\sqrt{|\alpha|}}\mathrm{i}$, respectively.
    \item[ \normalfont{Case }$\alpha>0$:] Let be ${\lambda}_{a+b}$ and ${\lambda}_{a-b}$ ($\mathbf{u}_{a+b}$ and $\mathbf{u}_{a-b}$) the eigenvalues (eigenvectors) of  the real matrices $\mathbf{A}+\dot{\mathbf{B}}$ and $\mathbf{A}-\dot{\mathbf{B}}$, respectively. 
    Then, the eigenvalues and eigenvectors of $\mathbf{X}$ are given by $\lambda=\lambda_1+\lambda_2\mathrm{i}$ and $\mathbf{u}=\mathbf{u}_1+\mathbf{u}_2\mathrm{i}$, respectively, with $\lambda_1=\frac{{\lambda}_{a+b}+{\lambda}_{a-b}}{2}$, $\lambda_2=\frac{{\lambda}_{a+b}-\lambda_1}{\sqrt{\alpha}}$, $\mathbf{u}_1=\frac{\mathbf{u}_{a+b}+\mathbf{u}_{a-b}}{2}$ and  $\mathbf{u}_2=\frac{\mathbf{u}_{a+b}-\mathbf{u}_1}{\sqrt{\alpha}}$.
\end{description}
\end{lemma}

\begin{remark}
If $\alpha<0$,  then the eigenvalues $\lambda$ and eigenvectors $\mathbf{u}$ are $(\alpha\beta)$-tessarines numbers,  where the components $\mathrm{j}$ and $\mathrm{k}$ vanish.  However, if $\alpha>0$,  they take  the form $\lambda=a_1+b_1\epsilon+c_1 \mathrm{i}+d_1\epsilon\mathrm{i}$ and $\mathbf{u}=\mathbf{a}_2+\mathbf{b}_2\epsilon+\mathbf{c}_2 \mathrm{i}+\mathbf{d}_2\epsilon\mathrm{i}$, where $\epsilon\mathrm{i}=\mathrm{i}\epsilon$. 
\end{remark}

 \subsection{SVD of an $(\alpha)$-Complex Matrix}
The following lemma establishes the SVD of an  $(\alpha)$-complex matrix, based on the SVD of a related complex matrix.

\begin{lemma}\label{SVDlema}
 Consider that $\mathbf{X}\in \mathbb{C}_{\alpha\beta}^{p\times q}$.  
 The following properties hold:
\begin{description}
\item[ \normalfont{Case }$\alpha<0$:] Let $\dot{\mathbf{X}}=\dot{\mathbf{U}}\mathbf{\Sigma}\dot{\mathbf{V}}^{*}$ be the SVD of the complex matrix  $\dot{\mathbf{X}}$ with $\dot{\mathbf{U}}={\mathbf{U}}_1+{\mathbf{U}}_2\epsilon$ and $\dot{\mathbf{V}}={\mathbf{V}}_1+{\mathbf{V}}_2\epsilon$.
    Then,  $\mathbf{X}=\mathbf{U}\mathbf{\Sigma}\mathbf{V}^{\texttt{H}_{-1}^1}$ is the SVD of $\mathbf{X}$ with matrices
    $\displaystyle{\mathbf{U}={\mathbf{U}}_1+\frac{{\mathbf{U}}_2}{\sqrt{|\alpha|}}\mathrm{i}}$ and    $\displaystyle{\mathbf{V}={\mathbf{V}}_1+\frac{{\mathbf{V}}_2}{\sqrt{|\alpha|}}\mathrm{i}}$.
    \item[ \normalfont{Case }$\alpha>0$:] Suppose that $q\leq p$.  Let $\mathbf{E}=\mathbf{U}_1\mathbf{\Sigma}_1\mathbf{V}_1^{\texttt{T}}$ and $\mathbf{F}=\mathbf{U}_2\mathbf{\Sigma}_2\mathbf{V}_2^{\texttt{T}}$ be the SVDs of the real matrices  $\mathbf{E}=\mathbf{A}+\dot{\mathbf{B}} $ and $\mathbf{F}=\mathbf{A}-\dot{\mathbf{B}} $, respectively.  The corresponding sets of singular values  are denoted by $\sigma_{1,1}, \sigma_{1,2},\ldots, \sigma_{1,r_1}$ and $\sigma_{2,1}, \sigma_{2,2},\ldots, \sigma_{2,r_2}$, respectively.  Moreover,  $i$th column of the matrices
        $\mathbf{U}_n$, $\mathbf{V}_n$,  for $n=1,2$,  are denoted by $\mathbf{u}_n(i)$ and $\mathbf{v}_n(i)$, respectively.

       Given two permutations $\{\tau_n(\sigma_{ni})\}_{i=1}^{r_n}$,  for $n=1,2$, the following matrices are defined:
        \begin{equation*}
        \begin{split}
        \dot{\mathbf{\Sigma}}_{n}&=\left[
          \begin{array}{cccc}
            \tau_n(\sigma_{n,1}) & 0 & \ldots & 0 \\
          0 & \tau_n(\sigma_{n,2}) & \ldots & 0 \\
            \vdots&   \vdots&   \vdots&   \vdots\\
0 &0 & \ldots & \tau_n(\sigma_{n,r_n})\\
            0 &0 & \ldots & 0\\
            \vdots&   \vdots&   \vdots&   \vdots\\
             0 &0 & \ldots & 0\\
          \end{array}
        \right] \\
        \dot{\mathbf{U}}_{n}&=[\mathbf{u}_{n}(\tau_n(\sigma_{n1})),\ldots,\mathbf{u}_{n}(\tau_n(\sigma_{nr_n})), \mathbf{u}_{n}(r_n+1),
    \ldots,\mathbf{u}_{n}(p)] \\
        \dot{\mathbf{V}}_{n}&=[\mathbf{v}_{n}(\tau_n(\sigma_{n1})),\ldots,\mathbf{v}_{n}(\tau_n(\sigma_{nr_n})), \mathbf{v}_{n}(r_n+1),
    \ldots,\mathbf{v}_{n}(p)]
    \end{split}
        \end{equation*}
  for $n=1,2$.
        Then,  the SVD of  $\mathbf{X}$  is given by $\mathbf{X}=\mathbf{U}\mathbf{\Sigma}\mathbf{V}^{\texttt{H}_{1}^1}$,   where
        \begin{equation}\label{svdcp}
        \begin{split}
        \mathbf{U}&=\frac{\dot{\mathbf{U}}_{1}+\dot{\mathbf{U}}_{2}}{2}
        +\frac{\dot{\mathbf{U}}_{1}-\dot{\mathbf{U}}_{2}}{2\sqrt{\alpha}}\mathrm{i}\\
          \mathbf{\Sigma}&=\frac{\dot{\mathbf{\Sigma}}_{1}+\dot{\mathbf{\Sigma}}_{2}}{2}
        +\frac{\dot{\mathbf{\Sigma}}_{1}-\dot{\mathbf{\Sigma}}_{2}}{2\sqrt{\alpha}}\mathrm{i}\\
          \mathbf{V}&=\frac{\dot{\mathbf{V}}_{1}+\dot{\mathbf{V}}_{2}}{2}
        +\frac{\dot{\mathbf{V}}_{1}-\dot{\mathbf{V}}_{2}}{2\sqrt{\alpha}}\mathrm{i}
        \end{split}
        \end{equation}
\end{description}
\end{lemma}

\begin{remark}
If $\alpha<0$,  the singular values, i.e,  the nonzero elements of $\mathbf{\Sigma}$ are positive real numbers.  However, if $\alpha>0$,  the eigenvalues  are $(\alpha\beta)$-tessarine numbers, where the  $\mathrm{j}$ and $\mathrm{k}$ components vanish.  Moreover,  in  case $p<q $,  a SVD for $\mathbf{X}$ can be obtained in a manner similar  to \eqref{svdcp}.
\end{remark}



\nocite{*}
\bibliographystyle{elsarticle-num-names}

\end{document}